\documentclass{amsart}[10pt,draft]
\usepackage[colorlinks=true,          
            linkcolor=blue,
            citecolor=red,
            urlcolor=blue]{hyperref}
\usepackage{amssymb,amsfonts,amsmath,hyperref,mathrsfs,latexsym,stmaryrd,graphicx}
\DeclareMathAlphabet{\mathpzc}{OT1}{pzc}{m}{it}

 \newtheorem{thm}{Theorem}[section]
 \newtheorem{Lemma}{Lemma}[section]
 
 \newtheorem{Prop}{Proposition}[section]
 \newtheorem{Ex}{Example}[section]
 \newtheorem{Cor}{Corollary}[section]

\usepackage {amsfonts,amssymb}
\usepackage{latexsym}
\usepackage{mathrsfs}
\input xy
\xyoption{all}

\newcommand {\Def}{\textrm{Def}}
\newcommand {\MC} {\textrm{MC}}
\newcommand {\Art}{\textrm{Art}_\CC}

\newcommand{\fart}{\textrm{FArt}_\CC}
\newcommand{\fun}{\textrm{Fun}}
\newcommand{\sets}{\textrm{Sets}}

\newcommand {\CC}{\mathbb{C}}
\newcommand {\FF}{\mathbb{F}}
\newcommand {\KK}{\mathbb{K}}
\newcommand{\NN}{\mathbb{N}}

\newcommand {\RR}{\mathbb{R}}
\newcommand {\SSS}{\mathbb{S}}

\newcommand{\YY}{\mathbb{Y}}
\newcommand{\ZZ}{\mathbb{Z}}

\newcommand {\bh}{{\bf h}}
\newcommand {\bk}{{\bf k}}

\newcommand {\bs}{{\bf s}}

\newcommand {\bF}{{\bf F}}
\newcommand {\bG}{{\bf G}}
\newcommand {\bH}{{\bf H}}

\newcommand {\bP}{{\bf P}}
\newcommand {\bQ}{{\bf Q}}

\newcommand{\cA}{\mathcal{A}}
\newcommand{\cB}{\mathcal{B}}
\newcommand{\cC}{\mathcal{C}}

\newcommand{\cF}{\mathcal{F}}

\newcommand{\cH}{\mathcal{H}}

\newcommand{\cK}{\mathcal{K}}

\newcommand {\cO}{\mathcal{O}}

\newcommand{\cQ}{\mathcal{Q}}
\newcommand{\cW}{\mathcal{W}}

\newcommand{\cih}{\mathpzc{h}}

\newcommand{\cI}{\mathcal{I}}

\newcommand{\cS}{\mathcal{S}}

\newcommand{\cU}{\mathcal{U}}

\newcommand{\cY}{\mathcal{Y}}

\newcommand{\scZ}{\scZ}

\newcommand{\fh}{\mathfrak{h}}

\newcommand{\fg}{\mathfrak{g}}

\newcommand{\fl}{\mathfrak{l}}
\newcommand{\fm}{\mathfrak{m}}

\newcommand{\fs}{\mathfrak{s}}
\newcommand{\fsu}{\mathfrak{su}}
\newcommand{\ft}{\mathfrak{t}}
\newcommand{\slt}{\mathfrak{sl}(2,\CC)}

\newcommand{\fz}{\mathfrak{z}}

\newcommand{\faff}{\mathfrak{aff}}

\newcommand {\dbar}{\overline{\partial}}
\newcommand {\zbar}{\overline{z}}
\newcommand {\dzbar}{d\overline{z}}

\newcommand {\shom}{\textrm{\underline{Hom}}}
\newcommand{\mhom}{\textrm{Hom}}
\newcommand {\send}{\underline{End} }
\newcommand {\mend}{\textrm{End}}

\newcommand {\ad}{\textrm{ad} }

\newcommand{\lspan}{\textrm{span}}
\newcommand{\img}{\textrm{Im }}

\newcommand {\cok}{\textrm{coker}}

\newcommand{\tildel}{\widetilde{\delta}}
\newcommand{\ctimes}{\otimes_\CC}
\newcommand{\sotimes}{\otimes_{\cO_X}}
\newcommand{\pic}{\textrm{Pic}}

\newcommand {\io}{\iota}
\newcommand {\fii}{\varphi}
\newcommand{\hookr}{\hookrightarrow}

\newcommand {\dd}{\textrm{d}}
\include{biblio}
\title{A Universal Family of Deformations for the Uniformising Higgs bundle }
\author{Peter Dalakov}
\address{Sector for Mathematical and Theoretical Physics, INRNE, Bulgarian Academy of Sciences, Tsarigradsko Chaussee 72, Sofia 1784, Bulgaria}
\date{\today} 
\newtheorem*{thma}{Theorem A}
\newtheorem*{thmb}{Theorem B}
\newtheorem*{thmc}{Theorem C}
\newtheorem{proposition}[thm]{Proposition}

\newtheorem{corollary}[thm]{Corollary}
\theoremstyle{remark}
\newtheorem{remark}[thm]{Remark}
\theoremstyle{observation}

\theoremstyle{definition}

\begin{document}

\begin{abstract}
Fix a simple complex Lie group $G$ and a principal $\mathfrak{sl}(2,\mathbb{C})$ subalgebra of $\textrm{Lie }(G)$.
Then the moduli space of semi-stable, topologically trivial $G$-Higgs bundles on a hyperbolic, spin Riemann surface
 acquires  a marked point.
This is   the unique $\mathbb{C}^\times$-fixed point on the Hitchin section.
We describe a universal analytic family of deformations  which  provides holomorphic Darboux coordinates
in a neighbourhood of the section.
 This is a special case of a more general deformation-theoretic
construction in the spirit of Kuranishi theory. As a toy example  of the latter we consider the
tautological family of centralisers over the Kostant slice.
\end{abstract}
\maketitle 
\setcounter{tocdepth}{1}
\tableofcontents 
  \section{Introduction}\label{introduction}
      \subsection{Motivating example}
The affine line of companion  matrices
\[
\Sigma =
\left\{ 
\begin{pmatrix}
 0&\alpha\\ 1&0\\
\end{pmatrix}, \alpha\in\CC\right\}\subset \slt 
\]
provides a section for  $-\det: \slt\to \CC\simeq \slt\sslash SL(2,\CC)\simeq \ft/(\ZZ/2)$, where
$\ft$  is the Cartan subalgebra of $\slt$.
Nigel Hitchin observed  (\cite{hitchin_sd}) that   $\Sigma$ can be promoted to  a $(3g_X-3)$-dimensional family of Higgs fields on the vector bundle
$K^{1/2}_X\oplus K^{-1/2}_X$, where $X$ is a  Riemann surface of genus $g_X\geq 2$. The Higgs fields in  this family are given by the above  formula but
 with
 $\alpha\in H^0(X,K^2_X)$.  This observation has  numerous far-reaching consequences and generalisations.
On the other hand,  the tautological family of   centralisers over $\Sigma$ is isomorphic to $T_{\ft/(\ZZ/2)}$
and 
 can be trivialised by
\[
\CC^2\simeq 
  \left\{ 
\begin{pmatrix}
 0&\alpha\\ 1&0\\
\end{pmatrix}, 
		  \begin{pmatrix}
		   0&\alpha\xi\\
		  \xi&0\\
		  \end{pmatrix}\right\} \subset \slt\times\slt .
\]
The trace gives a
    complex symplectic form
on $\slt\times \slt$,
and the above trivialisation provides Darboux coordinates for the bundle of centralisers.
Among other things, in this note we show how to 
construct a 
 $(6g_X-6)$-dimensional family of Higgs bundles  by twisting appropriately  the above formula.
      \subsection{Background}
Let $G$ be a simple complex Lie group, and  $X$  a smooth, compact Riemann surface
of genus at least two. A $G$-Higgs bundle on $X$ is a pair $(\bP,\theta)$, where $\bP$ is a holomorphic principal
$G$-bundle, and $\theta\in H^0(X,\ad\bP\otimes K_X)$.
The moduli space $M_{Dol}(G)$ of topologically trivial, semi-stable $G$--Higgs bundles on $X$ 
was constructed by Hitchin (\cite{hitchin_sd}, \cite{hitchin_sb})
 and Simpson (\cite{hbls},\cite{moduli2}). It
admits a proper map, $\chi$, called \emph{the Hitchin map},
  to a vector space, $\cB_\fg$, \hyperref[hitchin_base]{the Hitchin base}.
The Hitchin map admits a section which is a ``global analogue''
of the well-known ``Kostant slice'' from Lie theory. The latter generalises the notion of ``companion matrices'' and is a section of the adjoint quotient morphism
$\textrm{Lie } G =\fg\to \fg\sslash G$.
The Kostant section is not canonical, but depends on
a choice of  Lie-algebraic data (\hyperref[principal]{a principal $\slt$-subalgebra}). 
Similarly, Hitchin's section  depends on such a choice, as well as on  a choice of a theta-characteristic $\zeta= K_X^{1/2}$.
 We assume that these choices are fixed  once and for all, and hence shall talk about \emph{the} Hitchin section. 
The interested reader can find more details  in the original paper \cite{hitchin_teich}, as well as in
\cite{don-pan} or \cite{ngo}.
It should be noted that $M_{Dol}(G)$ is a holomorphic symplectic variety, $\chi$  is a complex Lagrangian fibration, 
and the  section is Lagrangian.

With the choice of a principal $\slt$ and a theta-characteristic the moduli space  acquires 
a marked point as follows.
 There is a natural $\CC^\times$-action on $M_{Dol}(G)$, given by $\lambda\cdot [(\bP,\theta)]=[(\bP,\lambda\theta)]$ (\cite{hbls}, \cite{hitchin_sd}).  
The marked point is
the unique $\CC^\times$-fixed point lying  on the (image of the) Hitchin section.
In the special case $G=SL(2,\CC)$ it appeared in Hitchin's original paper \cite{hitchin_sd} (Example 1.5 on p.8) as the first nontrivial  example of a stable Higgs bundle. 
 We call it ``the uniformising Higgs bundle'' (the terminology goes back to \cite{simpson_uniformisation}) since the Hermitian-Yang-Mills metric on this bundle is 
obtained from the uniformising metric of the curve $X$. In physics (the logarithm of) this metric is known as a
Toda field.
 This   is also an example of a ``system of Hodge bundles'' in the terminology
of \cite{hbls}. It  corresponds, by the non-abelian Hodge theorem (\cite{hbls})
 to a variation of Hodge strutures, and in fact,  a very special one, a $G$-oper.
The uniformising Higgs bundle carries a (regular) nilpotent Higgs field, i.e., belongs to the ``global nilpotent
cone'' $\chi^{-1}(0)$, the $0$-fibre of the Hitchin map. 
More about  systems of Hodge bundles and  uniformisation can be found  in \cite{hitchin_sd}, \cite{simpson_uniformisation}, \cite{hitchin_teich}, \cite{hbls},
\cite{simpson_iterated}. 

All these special properties of the uniformising Higgs bundle impose restrictions on its deformation theory. Their r\^ole is discussed  in
Section \ref{prelim}.

	\subsection{Results and contents of the paper}
The main result in this paper is contained in Section \ref{main_example}, where  we describe a universal analytic family of deformations of the uniformising Higgs bundle, with base the germ 
$\left( \cB_\fg\times\cB_\fg^\vee, 0\right)$. Here $\cB_\fg^\vee$
is the  dual vector space to the Hitchin base. Our family has the property that the holomorphic symplectic form on $M_{Dol}(G)^{reg}$
 induces   the canonical symplectic form on $\cB_\fg\times\cB_\fg^\vee$, so we
 obtain   holomorphic Darboux coordinates in an analytic neighbourhood 
of $[(\bP,\theta)]$.
 As a by-product, we obtain a formula for the flow of the Hitchin section under linear Hamiltonian functions on
$\cB_\fg$. This generalises an unpublished observation of C.Teleman for the case of structure group $GL(n,\CC)$ (\cite{teleman_langlands}).

Our approach uses several analytical pieces of data. First, we work with
  (analytic) differential graded Lie algebras (dgla), so holomorphic bundles are described in terms of
their Dolbeault operators. And second, we use a small amount of Hodge theory.
On the other hand, our final formulae  are polynomial  and are of Lie-algebraic origin, so a 
  purely algebraic description of the flow may also be feasible.

Section \ref{centralisers} is devoted to a   toy-version of the main  example: we give there a trivialisation of  the tautological family of centralisers
over  the Kostant slice.

In Section \ref{prelim} we make some general remarks about deformation theory via dgla's. We also describe the special features of
the controlling dgla and sketch a general strategy that one can follow in order to understand such deformation problems.

The results   from Sections \ref{centralisers} and \ref{main_example} are a consequence of the special form of
the dgla's controlling the corresponding deformation problems.
In Section \ref{symplectic_kur} we give sufficient conditions on the controlling dgla under which  similar (weaker) results hold.

The remaining sections are supplementary.   In Section \ref{Lie}
we recall  results from Lie theory and set up notation, and in Appendix \ref{Kuranishi_Theory} we review for our reader's convenience the basics of Kuranishi theory.
In Section \ref{conventions} we give a glossary of notation.

Our main results are as follows.

Let  $\MC_{L^\bullet}$ (respectively $\Def_{L^\bullet}$)  denote the  Maurer-Cartan (respectively,  deformation) functor
of a dgla $L^\bullet$,
 and let
$pr$ be the natural projection $\MC_{L^\bullet}\to\Def_{L^\bullet}$.  
 Suppose  $L^1=L'\oplus L''$
satisfies assumptions $(1)$, $(2)$, $(3)$ from Section \ref{symplectic_kur}.
The two inclusions (resp. projections) are denoted by $\io'$, $\io''$ (resp. $\pi'$, $\pi''$).
 Let
 $\cH^1=\cH'\oplus \cH''\subset L^1$ be ``harmonic representatives''  of $H^1(L^\bullet)$ (see \ref{Kuranishi_Theory}) and let  $\bH:L^1\to \cH^1$ be the corresponding  projection,
$\bH=\bH'+\bH''$.
Let $\FF_{L^\bullet}:\Art\to \sets$ be the formal Kuranishi map.
 We define a functor 
$\SSS_L= \MC_{L^\bullet}\cap \ker\left[(1-\bH')\pi' \right]:\Art\to\sets $.

\begin{thma}[\ref{symplectic_slice},\ref{triv_symp}]
Let $L^\bullet$ be a dgla with $L^3=0$, $H^2(L^\bullet)=0$ and  $L^1=L'\oplus L''$,  satisfying   $(1)$, $(2)$, $(3)$ from Section \ref{symplectic_kur}.
Let $P\io''$ be a splitting of $d_1'$ and $\pi:L^2\to \img d_1'$ a projection.
Assume that the formal series $\Gamma\in  L^1\widehat{\otimes}\varprojlim \textrm{Sym}^\bullet(\cH^{1\vee})/\fm^k $ defined by
$\Gamma(h,v):=(h, (1+P\pi\ad_h)^{-1}(v))$
 satisfies $[\Gamma,\Gamma]\in \img d_1'\widehat{\otimes}\varprojlim \textrm{Sym}^\bullet(\cH^{1\vee})/\fm^k$.
Then:
      \begin{itemize}
    \item The natural transformation
\[
\Phi: \SSS_L\to \underline{\cH}^1= \underline{\cH}'\oplus\underline{\cH}''
 \]
\[
 \Phi_A(h,v)= (h,(1+P\pi\ad_h)v)\in \cH^1\otimes \fm_A
\]
is an isomorphism in $\fart$ and $\Phi^{-1}=\Gamma$. The composition
 \[
\textrm{pr }\circ \FF_L^{-1}\circ \Phi \colon \SSS_L\to \Def_L\]
 is \'etale.
If moreover $L^\bullet$ is  normed and $\img d_1'\subset L^2$ is closed,  then $\SSS_L$ is prorepresented
by the germ $(\cS,0)$, where 
\[
 \cS = \MC(L)\cap \ker\left[(1-\bH')\pi' \right]\subset \cH'\oplus L'',
\]
and $\Phi: (\cS,0)\simeq (\cH^1,0)$.
    \item Suppose that $L'$ and $L''$ are  in (weak) duality by a pairing $\langle\ ,\ \rangle$ and let $\omega_{can}$ be the canonical symplectic form on $L^1$.
 Then $\Gamma^*\omega_{can} = \omega_{can}$, provided
 $\img P\pi\ad_h\subset \cH'^\perp$  for all $h\in\cH'$.
In the normed case, $\Phi:\cS\to \cH^1$ gives holomorphic Darboux coordinates on $(\cS,0)$.
\end{itemize}
\end{thma}

\begin{thmb}[\ref{major_lie}]
    Let $\fg$ be a simple complex Lie algebra,  $\{y,\fh, x\}$  a principal $\slt$ subalgebra,  and $P$  the canonical splitting of $\ad_y$ determined by it.
Let $\pi$ be the projection onto $\img\ad_y$, $\Sigma= y+\fz(x)$  the Kostant slice and $I$ the tautological family of centralisers. Then
\[
 \Phi : \cS \equiv I\vert_{\Sigma}\to \fz(x)\times \fz(y)
\]
\[
 \Phi(h,u)= (h,(1+P\pi\ad_h)u)
\]
is an isomorphism. Moreover,
\[
 \Phi(h,u)= (h, u +P[h,u])
\]
and
\[
 \Gamma(h,v):= \Phi^{-1}(h,v) = \left(h,\left(\sum_{k=0}^{\cih}(-1)^k\left(P\circ \ad_h \right)^k\right)(v)\right),
\]
where $\cih$ is the Coxeter number.  Finally, $\Gamma^*\omega_{can}=\omega_{can}$, where $\omega_{can}$ denotes
the canonical symplectic form on $\fg\times \fg$, as well as its restrictions to $I$ and $\fz(x)\times\fz(y)$.
\end{thmb}

\begin{thmc}[\ref{exponential}]
Let $(\bP,\theta)$ denote the uniformising Higgs bundle.
 The   notation and assumptions are from Section \ref{main_example}, in particular, we denote by $P$ the splitting of $\ad\theta$
induced by a  principal $\slt$-subalgebra.
 Consider the holomorphic family of Higgs bundles
\[
\Gamma:  \cH'\times \cH''  \longrightarrow  A^{1,0}(\ad\bP)\oplus A^{0,1}(\ad\bP),
\]
\[
 \Gamma(h,v)= \left(h,\Phi_h^{-1}(v)\right) = \left(h, \sum_{k=0}^{\cih}(-1)^k\left((s^{-1}P\ctimes 1)\circ ad_{h} \right)^k(v) \right),
\]
where $(h,v)\in \cH'\times\cH''\simeq H^1(L^\bullet)\simeq\cB_\fg\times\cB_\fg^\vee $
and
\[\Phi_h= 1+s^{-1}(P\ctimes 1)\pi\ad_h\in \mend (A^{0,1}(\ad\bP)).\]
The family $\Gamma$
is a miniversal deformation of the uniformising Higgs bundle $(\bP,\theta)$. An explicit description of $\cH'\times\cH''\subset A^1(\ad\bP)$ is given
in Theorem \ref{exponential}.

 There exists an open neighbourhood  $\cU\subset\cB_\fg\times\cB_\fg^\vee$
containing $0$,  
for which $\Gamma\vert_{\cU}$ is a universal deformation. Moreover, $\Gamma^\ast\omega_{can}= \omega_{can}$.
\end{thmc}
At this point it may seem utterly unclear why is it possible to describe such a family of deformations. In short, the reason is the very special nature
of our marked point, and, respectively, of the controlling dgla. In Section \ref{prelim}, after reviewing the basics of deformations via dglas,
we describe why our results are in fact natural.

\subsection{Acknowledgements} I would like to thank Tony Pantev, Carlos Simpson, Meng-Chwan Tan and Stefano Guerra for discussions, comments and feedback.

\vfill
\eject

\section{Deformations via dgla's}\label{prelim}
      \subsection{Basics}
We start with  some  remarks on  deformation theory via differential graded Lie algebras.
This  is by now very classical, and there are many great references.
The ones which seem both pedagogical and closest to our purposes are
  \cite{fukaya_def},
 \cite{goldman-millson}, \cite{maurer}, \cite{manetti_complex}.
 We state only the bare minimum of results and definitions,  without  
motivate them in any way. All vector spaces and tensor products are over $\CC$.

A \emph{differential graded Lie algebra} (dgla) is a triple $(L^\bullet,d,[\ , \ ])$. Here
$L^\bullet =\bigoplus_{k\in\NN} L^k[-k]$ is a graded vector space, endowed with a bracket
$[\ ,\ ]: L^i\times L^j\to L^{i+j}$. The bracket is graded skew-symmetric and
 satisfies a graded Jacobi identity. Finally, 
 $d:L^\bullet\to L^{\bullet+1}$ is a differential ($d^2=0$), which is a graded derivation of the bracket.
The \emph{set of Maurer-Cartan elements} in a dgla is the   zero set of the quadric $\cQ:L^1\to L^2$, $\cQ(u)=du+\frac{1}{2}[u,u]$.
We write  $\MC(L):= \cQ^{-1}(0)$.
To a dgla $L^\bullet$ one associates a Maurer-Cartan functor $\MC_{L^\bullet}: \Art\to \sets$, defined as
$$\MC_{L^\bullet}(A)= \MC(L^\bullet\otimes A) = \left\{ u\in L^1\otimes \fm_A :  du+\frac{1}{2}[u,u]=0 \right\}. $$
Given $\gamma\in L^1$, one can check (\cite{goldman-millson}, Section 1.3) that
$d_\gamma:= d+\ad \gamma\in Der^1L$ satisfies
$(d+\ad\gamma)^2=\ad\cQ(\gamma)$, and hence,  if $\gamma\in \MC(L)=\cQ^{-1}(0)$,
$d_\gamma$ is a differential, giving  a new dgla structure on $L^\bullet$.
 There is a Lie algebra homomorphism $L^0\to\faff(L^1)$ (the affine vector fields on
$L^1$), given by $\lambda\mapsto \left(\gamma\mapsto  - d_\gamma(\lambda)\right)$, and this affine vector field preserves the
set of Maurer-Cartan elements. We  define an
 action of $\exp(L^0\otimes \fm_A)$ on $L^1\otimes \fm_A$   by
\[  
 \exp(\lambda): u\mapsto \exp(\ad \lambda)(u)+\frac{I-\exp(\ad\lambda)}{\ad\lambda}(d\lambda)
\]
and define the deformation functor $\Def_{L^\bullet}: \Art\to \sets$ by 
\[
\Def_{L^\bullet} (A)= \MC_{L^\bullet} (A)/ \exp(L^0\otimes \fm_A) . 
\]
Then  $\MC_{L^\bullet}(A)$ can be considered as (the set of objects of) a groupoid, whose morphisms are determined by the gauge action;
this is often referred to as \emph{the Deligne-Goldman-Millson groupoid}.
Details about it  can be found in any of the references, e.g.,   Section 2.2. of \cite{goldman-millson}.
Deformation problems are described by  deformation functors $\Art\to \sets$, assigning to $A\in\Art$ the set of isomorphism classes
of deformations over $\textrm{Spec} A$. We say that a problem is
 governed (controlled) by a dgla, if its   deformation functor  is 
isomorphic to $\Def_{L^\bullet}$ for some dgla $L^\bullet$.

A dgla is called \emph{normed} (\cite{gm_kur}), if it is endowed with a norm, with respect to which $d$ and $[,]$
are continuous. It is called an \emph{analytic} dgla, if moreover it is endowed with continuous splitting $\delta$ ,  compatible 
with the other sructures.
See Appendix \ref{Kuranishi_Theory} or \cite{gm_kur} for the definition of splitting and details about compatibility.

If $L^\bullet$ is normed,
by a \emph{holomorphic family} of deformations of $\Def_{L^\bullet}(\CC)$ over a (pointed) complex manifold $(o,\cU)$ we mean a holomorphic map
$\Gamma: \cU\to \MC(L)\subset L^1$, $\Gamma(o)=0$. Holomorphicity  
makes sense even if $L^1$ is infinite dimensional, since it means continuous differentiability together with
$\CC$-linearity of $\dd \Gamma$. If $\cU$ is an open subset of a vector space and $\Gamma$ is a polynomial map, then
holomorphicity makes sense even if $L^1$ has no topology.

 One defines analogously deformations over a germ of an analytic subspace of $\CC^N$
or  more general analytic spaces,  see e.g. \cite{fukaya_def}, section 8.2. The Kodaira-Spencer map   
$\textrm{KS}: T_o \cU\to H^1(L^\bullet)$ is defined by $\textrm{KS}(\xi)=[\xi(\Gamma)(o)]$, where $\xi$ is thought of as a derivation.
    \subsection{Broader Context of the Paper}
One of the main outcomes of  this paper is the   explicit description of  a convenient universal family of deformations of a particular marked
point 
in a particular moduli space. Why is this possible at all?
The reason is the very special nature of our marked point, i.e., 
the very special form of the controlling differential graded Lie algebra
$L^\bullet$.
This ``speciality'' manifests itself in three ways:
\begin{enumerate}
 \item The dgla $L^\bullet$ is the total complex of a double complex. The bigrading (by Hodge type) and its interaction with the bracket put a restriction on
the set of Maurer-Cartan elements. 
 \item  The HYM  metric provides  a space
$L^1\supset \cH^1\simeq H^1(L^\bullet)$ of harmonic representatives. It comes with a  decomposition
$\cH^1\simeq \cH'\oplus \cH''$ into Lagrangian subspaces  which 
\emph{consist of  Maurer-Cartan elements}. The
 natural maps to $\Def_L$ are \'etale onto their images, see \ref{etale}. 
 \item The dgla $L^\bullet$ has an extra (finite length) grading  and one of the differentials of the double complex 
is   a shift with respect to it.
\end{enumerate}
Item $(1)$  holds for the dgla controlling the deformations of   any Higgs bundle.
Items $(2)$ and $(3)$ are related to the fact that $(\bP,\theta)$  is   a $\CC^\times$-fixed point and hence, by the non-abelian Hodge theorem,
 corresponds to a (polarised) $\CC$-VHS
(\cite{hbls}). The latter carries   two pieces of data: polarisation and Hodge filtration.   Item $(2)$ uses the particular
form of the polarisation and the Hodge structure on $H^1(L^\bullet)$. The grading in item $(3)$ is inherited from the  Hodge filtration on the associated
$\CC$-VHS.
It is well-known that for a smooth projective variety with a $\CC^\times$-action, the tangent space at an isolated fixed point decomposes into
``incoming'' and ``outgoing'' directions. The situation here is analogous, with  $\cH'\simeq \cB_\fg$ (resp. $\cH''\simeq \cB_\fg^\vee$) corresponding to  incoming
(resp. outgoing) directions.

We now argue that in such a situation there is a natural strategy for writing a (semi-universal) family of deformations.

For that we  look at the above three items
from a more general perspective, which
 is partially influenced  by
the  discussion of monads in  \cite{donaldson_kronheimer}  (Sections 3.1.3
and 3.2.1).

Suppose that $L^2$ and $L=L'\oplus L''$  are two complex vector spaces
and  that $\cQ:L\to L^2$ is an origin-preserving, ``off-diagonal'' quadratic map.
This means that
 $\cQ= \cQ_1 +\frac{1}{2}\cQ_2$, where $\cQ_1= \cQ_1'+\cQ_1''\in \mhom(L, L^2)$  and $\cQ_2\in\textrm{Hom}(L'\otimes L'', L^2)$.
Consider the quadric $M=\cQ^{-1}(0)\subset L$ and its ``tangent bundle'' $T_M=\ker \dd\cQ\vert_{M}\subset L\times L$,
$\dd\cQ_\lambda= \cQ_1 + \lambda\lrcorner\cQ_2$.
\footnote{If $L$ is infinite-dimensional and is not equipped with topology we do not have a notion of  vector bundle, hence the use of inverted commas.}
The quadric $M$ contains lots of affine spaces. In particular,  
$T_{M,0}=\ker\cQ_1$ contains two distinguished subspaces, $T'= \ker\cQ_1'\subset M$ and $T''=\ker\cQ_1''\subset M$.
The addition map $L\times L\to L$ identifies	
 $p_2^*L''\vert_{L'}\subset T_L$  with $L$
and
 $T_M\cap p_2^*L''\vert_{L'\times\{0\}}$ with 
 the ``slice''  $E:= M\cap (T'\times L'')$. The latter is a family of vector spaces parametrised by $T'$.
One may require (or look for conditions) that $E$
be a vector bundle (possibly, after suitable completion), so that all fibres  $E_h$, $h\in T'$  will be
 isomorphic to  $E_0=T''$.
 Suppose now $P:\textrm{Im}\cQ_1''\to L''$ is a splitting of the linear map $\cQ_1''$, and that 
$\pi:L^2\to \textrm{Im}\cQ_1''$ is a projection onto its image. Then the family of linear maps
$\Phi_h= 1 + P\pi (h\lrcorner\cQ_2)$ gives an identification $\Phi:E\simeq T'\times T''$, or rather, $E\vert_\cU\simeq \cU\times T''$, for some set $\cU$ around $0\in T'$,
determined by the condition that $\Phi_h$ be invertible. If $L$ is equipped with a topology in which the inverse function theorem holds, then $\cU$
can be taken to be an (analytic) open   set. If there is a (weak) duality pairing $L'\times L''\to \CC$, and $L$ is equipped with the corresponding canonical symplectic form $\omega_{can}$, 
then, under certain mild ``orthogonality'' conditions, $(\Phi^{-1})^\ast\omega_{can}=\omega_{can}\vert_{T'\times T''}$.

We shall apply this general strategy to the
 Maurer-Cartan quadric
  $\cQ(x)=dx+\frac{1}{2}[x,x]$. 
It is here that  Item $(3)$ enters: the choice of a principal $\slt$-subalgebra gives a natural splitting $P$ of 
the differential
$\cQ_1''=\ad\theta$, and the formal series for $\Phi^{-1}$ terminates due to the
finite length of the filtration, i.e., the
 nilpotence of $\theta$. Due to Item (2), we have harmonic representatives of $H^{1}(L^\bullet)$
and  can restrict $\Phi^{-1}$ to the subspace $\cH^1$.

This construction is formally similar to the standard construction of the Kuranishi family, but 
the r\^ole of 
 Green's operator $\bG$ is played by the much simpler splitting $P$.

\section{Symplectic Kuranishi Map}\label{symplectic_kur}

In this section we abstract some basic properties of the dgla's which occur in our examples of interest
and explore their deformation theory. 

Consider a dgla, $L^\bullet$,  whose $L^1$-term admits 
a non-trivial decomposition into a direct sum
$L^1=L'\oplus L''$, such that the two subspaces are:
\begin{enumerate}
 \item Isotropic for the bracket: $[L',L']=0=[L'',L'']$\label{condition1}
 \item Preserved under $\ad L^0$: $[L^0,L']\subset L'$, $[L^0,L'']\subset L''$.\label{condition2}
\end{enumerate}
Hence $(L^\bullet,d)$ contains as a subcomplex (not sub-dgla!)
the total complex of
\[
\xymatrix{L''\ar[r]^{d_1'} &L^2 \\ L^0\ar[u]^-{d_0''}\ar[r]^{d_0'}& L'\ar[u]^{d_1''}}. 
\]
If $L^\bullet$ is an analytic  dgla, we assume that the two subspaces $L'$ and $L''$ are closed, and the (co)product is in 
the category of topological vector spaces.
We denote by $d_k'$ the horizontal differentials, and by $d_k''$ the vertical ones. Notice that 
$\ker d_1'\subset L''$ and $\ker d_1''\subset L'$.
\begin{Ex}\label{lie_dgla}
Let $\fg$ be a complex Lie algebra and let
 $L^\bullet = \oplus_k L^k[-k]$, where  $L^0=\fg$, $L^1=\fg\oplus \fg$, $L^2=\fg$. Fix $y\in\fg$ and
endow $L^\bullet$ with   differentials $d_0=(\ad _y,0)^T$, $d_1=(0,\ad_y)$. There is a unique  bracket 
on $L^\bullet$ for which the above assumptions hold and which coincides with  the Lie bracket on  $L^0$.
\end{Ex}
\begin{Ex}\label{higgs_dgla}
 Let $X$ be a smooth compact curve, $E$ a holomorphic vector  bundle on it, and $\theta\in H^{0}(X,\send E\otimes K_X)$ a  Higgs field.
 Let $L^\bullet = \bigoplus_{p+q=\bullet} A^{p,q}(\send E)$ with differential $\dbar_E+\ad\theta$. Then we can take
$L' = A^{1,0}(\send E)$, $L''= A^{0,1}(\send E)$. The conditions on the bracket are satisfied for type reasons.
\end{Ex}
Finally, we  impose the following crucial   assumption
\begin{enumerate}
 \item[(3)]\label{condition3} Suppose $\img d_1'\subset L^2$ is split.
Suppose that $L^\bullet$ admits a splitting $\delta$ (see Appendix \ref{Kuranishi_Theory}) for which
  the direct sum decomposition of  $L^1$
 induces a non-trivial decomposition of $\cH^1$ into   $\cH^1 =  \cH'\oplus \cH''$ and $\cH''=\ker d_1'$. Fix one
such $\delta$. Denote by 
$\bH= \bH' + \bH'':L^1\to \cH^1$  the harmonic projection.
\end{enumerate}
\begin{Prop}\label{etale}
Suppose that  $H^2(L^\bullet)=0$ and $(1)$, $(2)$, 	$(3)$ hold.
Then there exist natural   morphisms in $\fart$
\[
\underline{\cH}'[-1] \subset\MC_{L^\bullet}\to \Def_{L^\bullet}
\]
\[
 \underline{\cH}''[-1] \subset\MC_{L^\bullet}\to \Def_{L^\bullet}
\]
which are \'etale onto their images.
\end{Prop}
\emph{Proof:}\\
By $(1)$ we have  
$\ker d'\subset \MC(L)$ and $\ker d''\subset \MC(L)$. But
$\cH'\subset\ker d''$ and $\cH''\subset\ker d'$, and the resulting inclusions
$\cH'\subset \MC(L)$ and $\cH''\subset\MC(L)$ induce the above-stated morphisms in $\fart$.
Since $H^2(L^\bullet)=0$ (obstructions vanish),
the Kuranishi map equals the identity on $\MC_L\cap \underline{\cH}^1 $, and the
(formal) Kuranishi functor $\KK_L$ equals $\underline{\cH}^1 $. 
 By  \cite{gm_kur}, Section 3 or \cite{maurer}, Theorem 4.7  (see also the Appendix \ref{Kuranishi_Theory}) we have an
\'etale morphism
\[\xymatrix@1{\KK_L\ar[r]^-{\FF^{-1}}& \YY_L=\MC_L\cap\ker \delta \ar[r]& \Def_L}. 
\]
\qed

We shall now digress and  make some   elementary remarks on   dgla's with vanishing  $L^3$ and $H^2(L^\bullet)$.
\begin{thm}
 Let $L^\bullet$ be a dgla with $L^3=0$ and $H^2(L^\bullet)=0$. Let $\widetilde{\delta}:L^2\to L^1$ be a splitting of $d_1$,
that is, $d_1\tildel=1_{L^2}$. Fix a subspace   $\cH^1\subset L^1$, isomorphic to $H^1(L^\bullet)$, and  
 consider  the formal power series
$\Gamma \in L^1\widehat{\otimes}\varprojlim \textrm{Sym}^\bullet(\cH^{1\vee})/\fm^k$,
$\Gamma= \sum _{k=1}^\infty \Gamma_k$, where
$\Gamma_k\in L^1\otimes \fm^k/\fm^{k-1}$ is defined inductively by
\begin{equation} \label{series2}
\Gamma_1(x) =x, \ 
 \Gamma_k(x)= -\frac{1}{2}\tildel\sum_{n=1}^{k-1}[\Gamma_n(x),\Gamma_{k-n}(x)].
\end{equation}
 Then $\Gamma$, thought of as a formal map $(\cH^1,0)\to L^1$,
 determines a formal miniversal family of deformations of $\Def_L(\CC)$ over $(\cH^1,0)$. 

If $L^\bullet$ is a normed dgla and  the above series converges in some  neighbourhood, $\cU$, 
of $0\in \cH^1$, then the corresponding family $\Gamma: \cU\to  L^1$ is a miniversal analytic family of deformations of $\Def_L(\CC)$.
\end{thm}
    \begin{remark}
In coordinates $\Gamma$ is described as follows.
   We fix  a  basis, $\{t_i\}$, $i=1...d$, of $(\cH^1)^\vee$.
 Then  
$\Gamma= \sum _{k=1}^\infty \Gamma_k \in L^1\widehat{\otimes} \CC\llbracket t_1,\ldots, t_d \rrbracket   $, and $\Gamma_k = \sum_{|J|=k}\Gamma_{J,k}t^J$, 
where $J$ is a multi-index.
    \end{remark}

\emph{Proof:} \\
This is a  statement about  power series which can  be related to some classical   deformation-theoretic  calculations  (see, e.g. \cite{kodaira_nirenberg_spencer} or 
\cite{kuranishi_complete}) . Since the proof is easy and instructive, we are going to give it here anyway.

 On one hand,
reading the Maurer-Cartan equation ``up to order $k$'' we see that any formal power series solution  has to satisfy
\[
 d_1\Gamma_k + \frac{1}{2}\sum_{n=1}^{k-1}[\Gamma_n,\Gamma_{k-n}]=0.
\]
On the other hand, applying $d_1$ to both sides of the proposed recursive formula for $\Gamma_k$ we get 
\[
 d_1\Gamma_k = -\frac{1}{2}d\tildel \sum_{n=1}^{k-1}[\Gamma_n,\Gamma_{k-n}] = -\frac{1}{2}\sum_{n=1}^{k-1}[\Gamma_n,\Gamma_{k-n}].
\]
By construction (\ref{series2}), the (formal or analytic) Kodaira-Spencer map of this family is the identity, 
and hence it is isomorphic to the Kuranishi family, which is a miniversal deformation. See e.g., \cite{fukaya_def} or
Appendix \ref{Kuranishi_Theory} for other references and comments. We emphasise  that this family \emph{need not}
be \emph{the} Kuranishi family.\\
\qed
\begin{remark}
 Kuranishi theory for dgla's (or $L_\infty$-algebras) is based on a choice of 
``splitting'' (or passing to a minimal model), see Appendix \ref{Kuranishi_Theory} for the relevant definitions. 
This involves a degree $-1$
endomorphism of $L^\bullet$, $\delta$, which in particular satisfies  $d\delta +\delta d = 1-\bH$, where $\bH$ is a ``harmonic projection''.
There is a well-known power-series solution of the Maurer-Cartan equation (the inverse of the formal Kuranishi map, see the Appendix \ref{Kuranishi_Theory}),
known from the works of Kuranishi, Kodaira-Nirenberg-Spencer, Huebschmann-Stasheff and many others. It 
 is given exactly by the  above formula (\ref{series2})
but
with $\delta$ instead of $\tildel$ (i.e., by (\ref{series})). The latter  formula  involves only $\delta_2$ and none of the
other $\delta_i$! 
To verify that the series (\ref{series})  gives a formal solution, one proceeds essentially as in the above proof. The main difference is that now
instead of $d_1\tildel =1$ we have
 $d_1\delta_2 = 1 -\delta_3 d_2  -\bH_2$.
But  $\bH_2=0$  since $H^2(L^\bullet)=0$, and the term involving $\delta_3 d_2$ vanishes due to the 
fact that $d$ is a derivation of the bracket, combined with associativity 
(graded Jacobi identity). If  $L^3=0$, this latter term is not present at all, so $d_1\delta_2 = 1$. Since the series
involves only $\delta_2$, we could start with any splitting (and ignore the remaining $\delta_i$) and will still get a formal
solution.
\end{remark}
We now return to our discussion of dgla's with a decomposition and state a version of the above theorem
based on splitting $d_1'$ only.

 Let $L^\bullet$ be a dgla with  $H^2(L^\bullet)=0$ and $L^3=0$,  satisfying the assumptions $(1)$,
$(2)$ and $(3)$.  
 Let $\pi: L^2\to \textrm{Im}d_1'$ be a projector and $P\io''$ a splitting of $d_1'$, so
 the linear map
$
 \tildel = \left(\begin{array}{c}
                 0\\ P\pi
                \end{array}
 \right) : L^2 \longrightarrow L^1$
satisfies $d_1\tildel = \pi$. 
\begin{thm} \label{triv_formal_analytic}
The formal power series
$\Gamma \in L^1\widehat{\otimes}\varprojlim \textrm{Sym}^\bullet(\cH^{1\vee})/\fm^k$ given by
\begin{equation}\label{series3}
 \Gamma (h,v) = \left( \begin{array}{r}
                        h\\  (1 + P\pi\ad_h)^{-1}(v)
                       \end{array}
 \right) = \left( \begin{array}{r}
                        h\\ v
                       \end{array}
 \right) + \sum_{k=1}^{\infty} (-1)^{k-1}\left( \begin{array}{r}
                        0\\ (P\pi\ad_h)^{k-1}(v)
                       \end{array}
 \right),
\end{equation}
$(h,v)\in \cH'\oplus \cH''$,
is a formal deformation of $\Def_{L^\bullet}(\CC)$ over $(\cH^1,0)$ if and only if  
$[\Gamma,\Gamma]\in \img d_1'\widehat{\otimes}\varprojlim \textrm{Sym}^\bullet(\cH^{1\vee})/\fm^k$. 

 If moreover $L^\bullet$ is a normed dgla and $\img d_1'\subset L^2$ is closed,  then 
there exists a neighbourhood of the origin,
 $\cU\subset \cH^1$,  such that
the family $\Gamma: \cU\to L^1$
is a miniversal analytic family of deformations of $\Def_{L^\bullet}(\CC)$.
\end{thm}
\emph{Proof:}\\
The formal statement is proved exactly  as in the previous theorem. Indeed, due to the isotropy of the bracket and
the choice of $\tildel$, the formula (\ref{series2}) reduces to the formula (\ref{series3}).
But since here $d\tildel=\pi$, we have  that $\Gamma$ satisfies
$d\Gamma = -\pi [\Gamma,\Gamma]$, which coincides with the Maurer-Cartan equation if and only if the right hand
side is $[\Gamma,\Gamma]$, i.e., $(\pi-1)[\Gamma,\Gamma]= 0$.

For the analytic statement notice that  the power series is essentially the geometric series, and since $\ad$, $\pi$ and $P$ are  continuous, 
the series will converge for $h$ sufficiently small, that is, for $x=(h,v)\in \cU= B_\epsilon\times \cH''$, where
$B_\epsilon\ni 0$ is  a ball of sufficiently small radius $\epsilon$. The Kodaira-Spencer map of the family is the identity, so it is miniversal by Theorem 1.3.3.,
\cite{fukaya_def}. 
\qed
\begin{corollary}\label{nilpotent_cor}
 Let $L^\bullet$ be as in the statement of  the  theorem, and 
assume that $\pi$, $P$, $\ad$ extend to continuous linear maps on some completion
 $\widehat{L}^\bullet$. Let $\delta$ be a compatibly chosen splitting of $\widehat{L}^\bullet$.
If $P\pi\ad_h$ is (locally) nilpotent for all $h\in \cH'$, then $\Gamma: \cH^1\to L^1\subset \widehat{L}^1$ is a miniversal analytic family of deformations of
$\Def_L(\CC)$.
\end{corollary}
\emph{Proof:}\\
If $P\pi\ad_h$ is locally nilpotent for all $h$,  then $(1+P\pi\ad_h)^{-1}(\cH'')\subset L\subset \widehat{L}$.\qed

\begin{proposition}\label{symplectic_slice}
Let $L^\bullet$ be a dgla satisfying  assumptions $(1)$, $(2)$, $(3)$.
Let $P\io''$ be a splitting of $d_1'$ and $\pi:L^2\to \img d_1'$ a projector.
Assume that  $L^3=0$, $H^2(L^\bullet)=0$, and  
$[\Gamma,\Gamma]\in \img d_1'\widehat{\otimes}\varprojlim \textrm{Sym}^\bullet(\cH^{1\vee})/\fm^k$. 

Let
$\SSS_L\in \fart$ be the functor  $\SSS_L= \MC_L\cap \ker\left[(1-\bH')\pi' \right]$. Then
\[
\Phi: \SSS_L\to \underline{\cH}^1= \underline{\cH}'\oplus\underline{\cH}''
 \]
\[
\Phi_A(h,v)= (h,(1+P\pi\ad_h)v)\in \cH^1\otimes \fm_A
\]
is an isomorphism in $\fart$ and $\Phi^{-1}=\Gamma$.
The composition 
\[
\textrm{pr }\circ \FF_L^{-1}\circ \Phi : \SSS_L\to \Def_L
\]
 is \'etale.
If additionally $L^\bullet$ is a normed and $\img d_1'\subset L^2$ is closed,  then $\SSS_L$ is prorepresented
by the germ $(\cS,0)$, where 
\[
 \cS = \MC(L)\cap \ker\left[(1-\bH')\pi' \right]\subset \cH'\oplus L'',
\]
and $\Phi: (\cS,0)\simeq (\cH^1,0)$.
\qed
\end{proposition}
\emph{Proof:}\\
We have that $(h,v)\in \MC_L(A)\cap \ker\left[(1-\bH')\pi' \right](A)$  $\Longleftrightarrow$
$h\in\cH'\otimes \fm_A$ and $v\in \ker (d'+\ad_h)$.

Now, $1+P\pi\ad_h$ maps
$\ker (d'+\ad_h)$ to $\ker d'$, since on  the former $\ad_h$ equals $-d'$, and $1+P\pi\ad_h$ equals 
$(1-Pd')$, the projector onto $\ker d'$.  
Since $\fm_A$ is nilpotent, $1+P\pi\ad_h$ is invertible for all $h$, and hence it maps injectively
$\ker (d'+\ad_h)$ to $\ker d'$. The condtion $(\pi-1)[\Gamma,\Gamma]=0$ means, by Theorem \ref{triv_formal_analytic}
 that
if $(h,v)\in \cH^1\otimes \fm_A= (\cH'\oplus\ker d')\otimes \fm_A$, then
$(h, (1+P\pi\ad_h)^{-1}v)\in \MC_L(A)$, i.e., belongs to $\SSS_L(A)$. Hence   $\Phi_A$	
is an isomorphism.

The composition $\textrm{pr }\circ \FF^{-1}\circ \Phi$ is \'etale  since $\Phi$ is an isomorphism and $\textrm{pr}\circ \FF^{-1}$ is \'etale by    \cite{gm_kur}, Section 3 or \cite{maurer}, Theorem 4.7.
\qed

For the next two corollaries,   assume that  
$\widehat{L}^\bullet$ is a normed dgla, which is the completion of a dgla
  $L^\bullet$ with respect to some norm. Also, assume that   $\delta$ is a (compatibly chosen) splitting  
 and that $P$, $\pi$ and $\ad$ extend to continuous operators on $\widehat{L}^\bullet$. 
\begin{corollary}
Let  both $L^\bullet$ and $\widehat{L}^\bullet$ satisfy the assumptions of the theorem.
Then, if $P\pi\ad_h$ is (locally) nilpotent for all $h\in \cH'$, the slice $\cS$ satisfies
$\cS = \Phi^{-1}(\cH^1)\subset L^1\subset \widehat{L}^1$.
\end{corollary}

\begin{corollary}
  Let $L^\bullet$  be a normed dgla satisfying the assumptions of the Theorem, except possibly the condition
$[\Gamma,\Gamma]\in \img d_1'\widehat{\otimes}\varprojlim \textrm{Sym}^\bullet(\cH^{1\vee})/\fm^k$.
Then $\Phi : (\cS,0)\simeq (\cH',0)\times \cH''$ is a holomorphic  vector bundle isomorphism over $(\cH',0)$ if and only if
 $[\Gamma,\Gamma]\in \img d_1'\widehat{\otimes}\varprojlim \textrm{Sym}^\bullet(\cH^{1\vee})/\fm^k$.
In particular, if $\dim L^\bullet <\infty$, $\dim\ker (d'+\ad_h)$ is constant in a (connected) neighbourhood of $0\in \cU\subset \cH'$ if and only if 
$\ad_h (1 + P\pi\ad_h)^{-1}(v)\in \textrm{Im }d_1'$ for all $h\in \cU$ and all $v\in \cH''$.
\end{corollary}

\emph{Proof:}\\
By the proof of  Theorem \ref{triv_formal_analytic}  $(1+P\pi\ad_h)$ is invertible for  $h\in B_\epsilon$, some $\epsilon>0$. By the inverse function theorem, its
inverse is analytic in some (possibly smaller) open, which we shall still denote by $B_\epsilon$.
 By the proof of Theorem \ref{symplectic_slice}  $\MC(L)\cap\ker\left[(1-\bH')\pi' \right]$
is a family of kernels, which we have trivialised by a holomorphic family of projectors. But by   \cite{banach_bundles}, \S 1, Theorem 1.5, Theorem 2.7 and  \S 3,
the image of a holomorphic family of projectors is a Banach vector bundle. Thus $\cS\cap B_\epsilon\oplus L''$ 
is  a Banach vector bundle precisely when $\Phi$ is an isomorphism.
Let us underline here that  $\dim H^\bullet(L^\bullet)<\infty$,  so both the base and the fibre of this vector bundle are finite dimensional vector spaces!
Of course, if  also $\dim L^\bullet<\infty$, then 
$\cS\cap B_\epsilon\oplus L''$ is a vector bundle if and only if $\textrm{rk}(1+P\pi\ad_h)=const$ on $B_\epsilon$.
\qed

 \begin{remark}
 Since $\cH^1$ is finite dimensional, probably  some clarification is needed reagrding the appearance of Banach vector bundles.
Our setup is  the following. We have a holomorphic  family  of linear maps between two (possibly) infinite dimensional vector spaces, $L''$ and $L^2$,
\emph{a priori} without topology: this is $\cH'\to \textrm{Hom}(L'', L^2)$,
$h\mapsto (d_1' + \ad_h)$.   We are interested in the collection of kernels, $\cS$.
We gave  conditions for the kernels to be of  finite, constant dimension ($=\dim \cH''$) and gave an explicit formal trivialisation, $\Phi$, of $\cS$.
If we want to put a topology on  $\cS$,  make it into  an honest vector bundle and have that  $\Phi$ be a  vector bundle trivialisation, then
we have to pass to a completion of $L^\bullet$. In the intended applications $\Phi_h$ is in fact a polynomial in $h$ due to nilpotence, and  
$\cS\subset L^1\subset \widehat{L}^1$.
 \end{remark}

We are ultimately interested in situations where
the ``local moduli space'' corresponding to $\Def_{L^\bullet}$ is symplectic, and the symplectic form
is induced by a (constant) symplectic form on $L^1$ for which the two subspaces $L'$ and $L''$ are isotropic.
The motivating example is the case when $L^1$ is a (weak) cotangent bundle. 

\begin{Lemma}
Let $L^\bullet$ be as in Theorem \ref{symplectic_slice}. Let $\omega$ be a skew-symmetric bilinear form on $L^1$
for which the subspaces 
$L'$ and $L''$ are isotropic. Then 
 $(\Phi^{-1})^*\omega$ 
 vanishes on the subbundles
$\cH^1\times\cH'\subset T_{\cH^1}$ and $\cH^1\times\cH''\subset T_{\cH^1}$. Moreover, $(\Phi^{-1})^*\omega$ on $\cH^1$ is invariant under translations along $\cH''$.
\end{Lemma}
\emph{Proof:} 
Let $\widetilde{\Phi}^{-1}$ denote the holomorphic map $\cH'\to \textrm{End }(L^2)$, $h\mapsto \Phi_h^{-1} = (1+P\pi\ad_h)^{-1}$.
Then $(\dd\Phi^{-1})_{(h,v)}(\xi',\xi'')= (\xi',0) + (0,(\dd\widetilde{\Phi}^{-1})_h(\xi'')(v))+(0,\Phi_h^{-1}(\xi''))$, so $\dd\Phi^{-1}$ preserves
$\ker d_1'$ and $\ker d_1''$ and the first statment  follows.  But the second of the three terms
vanishes identically due to  assumption $(1)$: $(\dd\widetilde{\Phi}^{-1})_h(\xi'')(v)=-P\pi[\xi'',v]=0$, and so
$(\Phi^{-1})^*\omega_{(h,v)}$ is independent of $v\in \cH''$.
  \qed
\begin{proposition}\label{triv_symp}
Let $L^\bullet$ (resp. $\widehat{L}^\bullet$) be a dgla, satisfying the assumptions of Theorem \ref{symplectic_slice}.
Suppose that $L'$ and $L''$ are placed in (weak) duality by a pairing $\langle\ ,\ \rangle$ and let $\omega_{can}$ be the canonical symplectic form on $L^1=L'\oplus L''$.
If $\textrm{Im }P\pi\ad_h \subset \cH'{^\perp}$ for all $h\in \cH'$,  then $(\Phi^{-1})^*\omega_{can} = \omega_{can}$.
In the normed case, $\Phi:\cS\to \cH^1$ gives holomorphic Darboux coordinates on $(\cS,0)$.
\end{proposition}
\emph{Proof:}\\ 
The canoncial symplectic form on $L^1$ is $\omega_{can}((\xi',\xi''),(\eta',\eta''))= \langle\xi',\eta'' \rangle -\langle\xi'',\eta'\rangle$.
Using the formula for $\dd \Phi^{-1}$ from the previous Lemma and the isotropy of $\cH'$ and $\cH''$ we get
\[
 (\Phi^{-1})^\ast\omega_{can}((\xi',\xi''),(\eta',\eta''))=\langle\xi', \Phi^{-1}(\eta'') \rangle -\langle  \Phi^{-1}(\xi''),\eta'\rangle.
\]
Substituting
$(1+P\pi\ad_h)^{-1} =1 -P\pi\ad_h(1+P\pi\ad_h)^{-1}$ into the previous formula and using the orthogonality assumption we get
$(\Phi^{-1})^*\omega_{can}((\xi',\xi''),(\eta',\eta''))= \langle\xi',\eta'' \rangle -\langle\xi'',\eta'\rangle$.
\qed

\begin{remark} 
With the above assumptions, $\cS$ is a Lagrangian foliation, with space of leaves (the germ of) $\cH'$.
Such a foliation carries a torsion-free flat connection along the leaves. Since  $T^\vee_{\cH'}\simeq \cS$
(as symplectic manifolds),  the affine structures on the leaves is induced by the vector space structure on the fibres,
and we have described it in terms of the controlling dgla.
\end{remark}
\section{Lie-algebraic preliminaries}\label{Lie}
We review here  some relevant facts from Lie theory mostly to set up notation. Details can be found  in
 \cite{cg} or \cite{kostant}.
Let $G$ be a simple complex Lie group, $\fg=\textrm{Lie}(G)$ and  $\textrm{rank}(\fg)=l$.
An element of $\fg$ is \emph{regular} if its centraliser is of  the smallest possible dimension, $l$.
 An element $\varphi\in \fg$ is \emph{semisimple} (respectively, \emph{nilpotent})
if $\ad_\varphi\in \mend(\fg)$ is semisimple (respectively, nilpotent). 
If  $\fg=\fs\fl(l+1)$,  the regular elements are trace-free matrices 
with   a single Jordan block per eigenvalue.  A regular nilpotent $\varphi$ is one which is conjugate to a single Jordan block with
zeros on the diagonal. We denote by $\fg^{reg}$, $\fg^{ss}$ and $\fg^{reg,ss}$ the sets of regular, semisimple and regular semisimple elements
of $\fg$.
One has
$\fg^{reg,ss}\subset \fg^{reg}\subset\fg$ and $\fg\backslash\fg^{reg,ss}\subset \fg$ is a divisor  while
$\fg\backslash\fg^{reg}\subset \fg$ is of codimension 3.

 The  notion of a regularity  makes sense for  reductive Lie algebras as well. In particular,  if   $\fii\in\fg\fl(n,\CC)^{reg}$, its centraliser  
$\fz(\varphi)$ is spanned by $\{\varphi,\varphi^2,\ldots, \varphi^{n}\}$. We  do not have such a convenient description of the centraliser for other Lie algebras.

By the Jacobson-Morozov lemma any nilpotent  $x\in \fg$   can be embedded in an $\slt$-subalgebra of $\fg$.
A \label{principal} \emph{principal $\slt$ subalgebra} is one which is spanned by
  two \emph{regular} nilpotent elements,   $x$ and $y$,
 and a semisimple $\fh\in\fg$.
The inclusion $\slt\hookr \fg$ exponentiates
 to a homomorphism $\varrho: SL(2,\CC)\to G$, a ``principal homomorphism''. The maximal compact $SU(2)\subset SL(2,\CC)$ maps to a
compact form of $G$.
Under the adjoint action of $\slt$  $\fg$ decomposes into
$l$ odd-dimensional irreducible representations:
\begin{equation}\label{decomp}
\fg = \bigoplus_{i=1}^{l} W_{m_i},\  W_{m_i}= Sym^{2m_i}(\CC^2),
\end{equation}
where $\CC^2$ is the standard representation of $\slt$. 
 The  spaces $W_{m_i}$ are $(2m_i+1)$-dimensional, so the restriction 
$SL(2,\CC)\to Aut(W_{m_i})$ of the adjoint representation to each $W_{m_i}$ factors through $PGL(2,\CC)$.
The restriction to the maximal compact makes $W_{m_i}$ into a representation of
$PSU(2)=SO(3)$. 
On each $W_{m_i}$ the eigenvalues of $\ad_{\fh}$ are even integers $2m$, where $-m_i\leq m \leq m_i$. The highest weight vectors span the
 centraliser  $\fz(x)$.   We shall
label the eigenspaces  by \emph{half} of the corresponding eigenvalue and shall
   let  $\fg_m$ stand for the eigenspace of  $\ad_{\fh}$ with  eigenvalue $2m$. The decomposition
\begin{equation}\label{lie_grading}
\fg = \bigoplus_{m=-\cih}^{\cih} \fg_m,  
\end{equation}
  is  called \emph{the principal grading} of $\fg$. 
The filtration $\cW_\bullet\fg$, $\cW_p\fg = \oplus_{2m\leq p}\fg_{m}$ is the 
 canonical (Deligne) filtration of the nilpotent endomorphism $\ad_y$. Intersecting
\ref{decomp} and \ref{lie_grading} we get 
a bigrading $\fg = \oplus \fg_{k,i}= \fg_k\cap W_{m_i}$. Then $\fz(x)=\oplus_i \fg_{m_i,i}$ and $\fz(y)=\oplus_i \fg_{-m_i,i}$.

The numbers $m_i$ are the \emph{exponents} of $\fg$ (or $G$).
For a simple Lie algebra they are all distinct except if $\fg=D_{2n}$, when the largest exponent has multiplicity two.
 We order the exponents, so that 
$m_i\leq m_j$ for $i<j$ and for the most part we shall write $W_i$ instead of $W_{m_i}$. In particular,
as $G$ is simple, $m_1=1$ and $W_1=\slt$ is the principal subalgebra.

The motivating  example is the $l$-th symmetric power embedding $\slt\hookr \fs\fl(l+1,\CC)$. Notice that it maps
 the standard generators $\{y_0 = E_{21}, h_0= E_{11}-E_{22}, x_0= E_{12}\}$ of $\slt$ to 
 the $(l+1)\times (l+1)$ matrices $\{y, h, x\}$, where $y = \sum_{p=1}^l E_{p+1,p}$, 
$h= \sum_{p=1}^{l+1}(l-2p+2)E_{p,p}$ and  $x= \sum_{p=1}^{l}p(l-p+1)  E_{p,p+1}$. In particular, $x\neq y^T$!

Let $\CC[\fg]^G\subset\CC[\fg]$ be  the ring of  $G$-invariants for the adjoint action.
The GIT quotient  is $\fg\sslash G:= \textrm{Spec } \CC[\fg]^G$, and
its points correspond to \emph{closures} of $G$-orbits. 
The closure of each $G$-orbit contains a unique open (regular) and a unique closed (semisimple) orbit.
By a theorem of Chevalley, $\CC[\fg]^G$ is isomorphic to a polynomial ring, i.e.,  $\fg\sslash G$ is non-cannonically isomorphic to 
a vector space.
 We can fix one such  isomorphism 
by choosing a basis for the $G$-invariant polynomials on $\fg$, say  $\{p_1,\ldots,p_l\}$,  $\deg(p_i)= m_i +1$.
We assume that our choice of invariant polynomials is
compatible with the decomposition $\fz(x)= \oplus_i \fz(x)\cap W_i$ induced by the 
 principal subalgebra. This means  that there exists a basis for $\fz(x)$ consisting of highest weight vectors 
$v_i\in W_i\cap \fg_{m_i}$,  such that $p_i(y+a_1v_1+ \ldots + a_lv_l)=a_i$.
This gives an identification $\CC[\fg]^G\simeq \CC[p_1,\ldots,p_l]$ and the Chevalley projection 
 $\fg\to \fg\sslash G$  can  be interpreted as a map $\fg\to \CC^l$. 
For $\fg=\fg\fl(n,\CC)$ this map sends a 
matrix to the (non-leading) coefficients of its characteristic polynomial.

Let  $\ft\ni \fh$ be a Cartan subalgebra and   $W$  the corresponding  Weyl group.
Chevalley proved that  $\ft\hookr \fg$  induces an isomorphism $\ft/W\simeq \fg\sslash G$.
 In \cite{kostant} it is shown that the adjoint quotient $\fg\to \ft/W$ becomes an isomorphism when restricted to \emph{the Kostant slice}
 $\Sigma= y + \fz(x)\subset \fg^{reg}$. Thus  $\Sigma$  provides a splitting $\ft/W\to \fg$ of the 
Chevalley projection.
 We shall also write $\bs$ for the affine-linear map $\bs: \fz(x)\to \Sigma$, $\bs(a)=a+y$.

 We shall use   one particular  principal $\slt$-subalgebra $\{y,\fh,x\}$ which 
 is  the standard one in the literature 
on opers (\cite{frenkel}) and which we describe now.

Fix  
 Chevalley generators
 $\{f_i, h_i, e_i\}$,  $\ft= \lspan \{h_i\}$, $i=1\ldots l$,
and assume  $\kappa(e_i,f_i)>0$. Fix positive roots $\Delta^+$.
Let  $\rho^\vee= \sum_i \rho_i^\vee h_i $ be the dual Weyl vector, i.e.,  half the sum of the positive coroots.
We take  $y = \sum_i f_i$, a regular nilpotent element, and $\fh=2\rho^\vee\in \ft$.
The   unique   $x$ for which  $\lspan\{x,2\rho^\vee,y\}\simeq \slt$ is 
     $x= \sum_i 2\rho^\vee_i e_i$.

The choice of Chevalley generators determines a split and a compact real  form of $\fg$ (\cite{bourbaki_lie_cpct},IX.16 \S 3). The former is  the real subalgebra   generated by
$\{e_i,f_i,h_i\}$. The latter is the $+1$ eigenspace of the \emph{anti-linear} extension, $\eta$,  of
$e_j\mapsto - f_j$, $f_j\mapsto - e_j$, $h_j\mapsto -h_j$. The   $+1$ eigenspace
 is generated by   $\{ih_j, e_j-f_j, i (e_j+f_j)\}$.
 Notice that $\eta$ is only a vector space involution, 
not a Lie algebra one. In the classical cases, $-\eta$ is  hermitian conjugation, so we often write $u^*$ for $-\eta(u)$.
 Our special choice of a  principal $\slt$ need  not  be 
``aligned''  with the choice of generators: the standard copy of $\fsu(2)\subset \slt$
 is \emph{not} mapped to the compact form   determined by the generators $\{e_i,h_i,f_i\}$.
All compact forms are conjugate, and  it is  easy to determine
the anti-linear involution preserving this one.
It  is the
anti-linear extension of
$e_j\mapsto - \frac{1}{2\rho^\vee_j}f_j$, $f_j\mapsto -  2\rho^\vee_j e_j$, $h_j\mapsto -h_j$. Its   $+1$ eigenspace
 is spanned  by   $\{\frac{i}{2\rho_j}h_j, e_j-\frac{1}{2\rho_j}f_j, i (e_j+\frac{1}{2\rho_j}f_j)\}$.

 The pairing  $(u,v) = -\kappa(u, \eta(v)) =\kappa(u,v^\ast)$ is an hermitian inner product on $\fg$, and
 $ad_u^* = ad_{u^*}$, where $\ad_u^\ast$ is the adjoint of $\ad_u\in \mend{\fg}$ with respect to it.
It is probably well-known that the different irreducible representations $W_i$ are orthogonal with respect to this inner product, but for lack
of reference we have  proved it in \cite{thesis}.

 Notice that by construction
  the principal $\slt$  (and all the representations $W_i$) are all real with respect to $\eta$ and in particular
$y^*=x,  \fh^*=\fh$. 
\section{Universal Centralisers}\label{centralisers}
Consider now   the tautological family of centralisers of regular elements
 \[
I=\left\{(v,u): [v,u]=0, v\in\fg^{reg},u\in \fg \right\}\subset \fg^{reg}\times\fg.
   \]
  The projection
 $\textrm{pr}_1:I\to \fg^{reg}$ makes this locally closed subvariety  into  a   rank $l$ vector bundle, a  subbundle of the trivial bundle $T_{\fg^{reg}} $.
The group $G$ acts on  $I$  diagonally by the adjoint action,  and the quotient 
is     \emph{the universal centraliser}. It is a hamiltonian  reduction of $T^\vee_{\fg}\simeq_\kappa T_\fg$,
and 
$I\sslash G\simeq T_{\ft/W}$ is a symplectic isomorphism. On the other hand, 
 $ I\sslash G\simeq I\vert_{y+\fz(x)}= \bs^*I$, and
we shall see  that the choice of a  principal  $\slt$ provides  a natural trivialisation 
$\bs^*I\simeq \fz(x)\times \fz(y)$, with the  the property that the symplectic form on $I\vert_{y+\fz(x)}\subset T_\fg$
pulls back to the standard symplectic form on $\fz(x)\times \fz(y)$.

The subspace $\fz(x)\simeq \textrm{coker}(\ad_y)$ provides a splitting, $P\in\mhom(\textrm{Im }\ad_y,\fg)$, of $\ad_y$.
To compute $P$ in examples  one can use that  each $W_{m_i}$ is an irreducible
$\slt$-representation, so a suitable multiple of $\ad_x$  inverts 
$\ad_y$ on $\img(\ad_y)$. 
For the actual coefficient, depending on $m_i$ and $k$, see  \cite{fulton_harris}, Lecture 11. 
The bigrading of $\fg$ provides  natural projections  $\pi: \fg\to \textrm{Im }\ad _y$ and 
$p^r_p:\fg_r\to \fg_{r,p}$. Note that $\pi, p^r_p\in \mend_0(\fg)$, while $P\in \mend_1(\fg)$.
Consequently,  for all  $h\in \fz(x)$, $P\pi\ad_h\in\mend_2 (\fg)$ and is hence nilpotent.
Note in passing that in this setup we also have a natural splitting of $\ad_x$, say $Q\in \textrm{Hom}(\img(\ad_x),\fg)$, $\ad_x\circ Q=1$.

We now formulate a technical Lemma.
\begin{Lemma}\label{induction}
  Let $\fg$ be a simple complex Lie algebra,  $\{y,\fh, x\}$  a principal $\slt$ subalgebra, and  $P$  the canonical splitting of $\ad_y$ determined by it.
Let $0\neq h\in\fz(x)$. Then, $\forall k\geq 0$, 
\[
 \ad_h(P\ad_h)^k(\fz(y))\subset \img \ad_y.
\]
Equivalently,  $\forall k\geq 0$, $(P\pi\ad_h)^k(\fz(y)) = (P\ad_h)^k(\fz(y))$.
\end{Lemma}

\emph{Proof:}\\
 For notational simplicity assume that $\fg\neq D_{2n}$. This is the only  simple Lie algebra  with a repeated exponent
(the largest exponent appears twice), and in that case the proof is exactly as the one that follows below, 
but one has to choose the two $W_i$'s corresponding to the maximal exponent in a way that they be orthogonal 
with respect to the inner product induced by the Killing form.

We work by   induction on $k$, and use
an observation from Clebsch-Gordan theory of $\textrm{SL}(2,\CC)$
(\cite{hitchin_teich}, p.458) regarding
 commutators of elements from different
$W_i$. Namely, 
$$pr^{m+n}_p\left( [\fg_{m,i},\fg_{n,j}] \right) = 0  \textrm{ unless } m_i + m_j + m_p = 1\textrm{ mod } 2 $$

For the base case $k=1$ we have to show that $\ad_h: \fz(y)\to Im(ad_y)$.
Let  $v\in\fz(y)$. Since $\fz= \oplus_i \fz(y)\cap  W_i$, and similarly for $\fz(x)$,
  we may assume $v\in \fg_{-m_j,j}\subset \fz(y)$ and
$h\in V_{m_i}=\fg_{m_i,i} \subset \fz(x)$.
Then
 $[h,v]=[e_{m_i},e_{-m_j}]\in \fg_{m_i-m_j}$, where $e_{m_i}$ (resp.  $e_{-m_j}$) is  a highest (resp.  lowest) weight vector in
$W_i$ (resp.  $W_j$). 
We claim that this commutator can never be in some $V_{m_p}=\fg_{m_p,p}$, that is
$pr^{m_i-m_j}_{m_i-m_j}([e_{m_i},e_{-m_j}])=0 .$
Indeed, if there were such a term,
 there would be an exponent $m_p$, such that
 $m_i-m_j = m_p$ and $m_i + m_j+m_p =1 \textrm{ mod }2$, which
would mean that $2m_i=1\textrm{mod }2$. So the base case is proved and  $P\circ\pi\circ\ad_h(v)= P\circ \ad_h(v)\in \fg_{m_i-m_j+1}$.

For the  inductive step, let  $(P\circ\pi\circ\ad_h)^k(v)=\left( P\circ ad_\eta \right)^k(v)$,  $k\geq 1$.
Then we can write it as a linear combination of elements in the $(km_{i}-m_j +k)$-th graded piece of $\fg$. 
 Such an element has  a nonzero projection  in some $W_{p_k}$ if
\[
m_{i}+m_j+m_{p_1}=2l_1+1,\  m_{i}+m_{p_1}+m_{p_2}=2l_2+1,\ldots m_{i}+m_{p_{k-1}}+m_{p_k}=2l_k+1,
\]
where $l_r\in\ZZ$.
 Adding  these up gives 
 \begin{equation}
 \label{parity1}
 km_{i} + m_j + 2 \sum_{s=1}^{k-1}m_{p_s} + m_{p_k} = \sum_{r=1}^k 2l_r+1
 \end{equation}
If  $\ad_h\left( P\circ ad_h \right)^k(v)$ has a nonzero projection  in some $V_{m_l}$, then it must be the case that
 $$(k+1)m_{i}- m_{j} + k = m_l,\ m_{i} + m_{p_k}+ m_l = 1 \textrm{ mod }2,$$
 and adding these we get
 \begin{equation}
 \label{parity2}
  (k+2)m_i  + m_{p_k} - m_{j} + k = 1\textrm{ mod }2.
 \end{equation}
 Finally, adding (\ref{parity1}) and (\ref{parity2}) we get an equality of the form
 $\textrm{even} + k =   \sum_{r=1}^{k+1} \textrm{odd}_r,$
 which is impossible since $k$ and $k+1$ are always of opposite parity.
\qed

We now reconsider    Example \ref{lie_dgla} with $y$ being
 the regular nilpotent element from the
principal $\slt$-subalgebra. 
Then 
$H^0(L^\bullet)=\fz(y)$, $\ker d_1=\fg\oplus \fz(y)$, $\img(d_0)= (\img(\ad_y),0)$ and so 
$H^1(L^\bullet)\simeq \fz(x)\oplus \fz(y)$. Thus $L^\bullet$ satisfies the assumptions $(1)$, 
$(2)$,$(3)$ from \ref{symplectic_kur}, with $\cH'=\fz(x)$ and $\cH''=\fz(y)$. We also have
\[
\MC(L) = \{(h,v): v\in \fz(y+h) \} = (\bs\times 1)^\ast I \subset \fg\times \fg
\]
\begin{thm}\label{major_lie}
   Let $\fg$ be a simple complex Lie algebra,  $\{y,\fh, x\}$  a principal $\slt$ subalgebra,  and $P$  the canonical splitting of $\ad_y$ determined by it.
Then
\[
 \Phi : \cS \equiv I\vert_{\Sigma}\to \fz(x)\times \fz(y)
\]
\[
 \Phi(h,u)= (h,(1+P\pi\ad_h)u)
\]
is an isomorphism. Moreover, 
\[
 \Phi(h,u)= (h, u +P[h,u])
\]
and
\[
 \Gamma(h,v):= \Phi^{-1}(h,v) = \left(h,\left(\sum_{k=0}^{\cih}(-1)^k\left(P\circ \ad_h \right)^k\right)(v)\right),
\]
where $\cih$ is the Coxeter number.  Finally, $\Gamma^*\omega_{can}=\omega_{can}$, where $\omega_{can}$ denotes
the canonical symplectic form on $\fg\times \fg$, as well as its restrictions to $I$ and $\fz(x)\times\fz(y)$.
\end{thm}

\emph{Proof:}\\
Since  $P\pi\ad_h$ is nilpotent  for all $h\in \fz(x)$, then $1+P\pi\ad_h$ is invertible.
Then $\Phi$ is an isomorphism by Theorem \ref{triv_formal_analytic} applied to   Example \ref{lie_dgla}: the series (\ref{series3})  is convergent
for all $h\in\cU=\fz(x)$ and  is actually a polynomial of degree at most $\cih$. But by Lemma \ref{induction} this polynomial
equals the one from the statement of the Theorem. The condition $(\pi-1)[\Gamma,\Gamma]=0$ clearly holds,
since all elements from the Kostant slice are regular.
Finally, the statement about the symplectic form holds by
Proposition \ref{triv_symp}, applied  to the dgla under consideration. Indeed, the Killing form is non-zero only on
$\fg_{m_i,i}\times \fg_{-m_i,i}$ and $\fg_{-m_i,i}\times \fg_{m_i,i}$, while   $P\pi\ad_h\in\mend_2(\fg)$, so the
orthogonality condition
 from \ref{triv_symp} is satisfied. \qed

\begin{Ex}
 Let $\fg=A_1=\slt$
with the standard generators $y=E_{21}$, $2\rho^\vee=E_{11}-E_{22}$, $x=E_{12}$. Then
$P(y)=\rho^\vee$, $P(h)=-x$ and
$\Gamma:\CC^2\simeq \fz(x)\times \fz(y)\to \slt\times\slt$ is given by
\[
 \Gamma(h,v)= (y+ \alpha x, \xi y+ \alpha\xi x)= \left(
\begin{pmatrix}
 0&\alpha\\ 1&0\\
\end{pmatrix}, 
		  \begin{pmatrix}
		   0&\alpha\xi\\
		  \xi&0\\
		  \end{pmatrix}
\right), \ \alpha,\xi\in\CC.
\]

\end{Ex}

\section{The uniformising Higgs bundle}\label{main_example}

\subsection{The Uniformising Higgs bundle}
Let us fix a theta-characteristic $K^{1/2}_X$. This is   a line bundle $\zeta\in\textrm{Pic}^{g-1}_X$, together with an isomorphism
$\zeta^{\otimes 2}\simeq K_X$. It is well-known that such a $\zeta$ always exists: $\zeta$ is a  spin-structure  
and $X$ is spin, since $w_2(X)=0$. There are $2^{2g}$ choices
of $\zeta$: the different theta-characteristics form a torsor over the points of order 2 in $\pic^0_X$.
Consider the $SL(2,\CC)$-Higgs pair   $(\zeta\oplus \zeta^{-1},\theta_0)$, 
 $\theta_0=\left( \begin{array}{rr}
                         0& 0\\ 1&0\\
                        \end{array}\right)$, where  $1$ is considered as a global section of $\zeta^{-2}\otimes_{\cO_X} K_X$.
Consider then  $\underline{Isom}(\zeta\oplus \zeta^{-1},\cO^{\oplus 2}_X)$, the $SL(2,\CC)$-frame bundle of $\zeta\oplus \zeta^{-1}$,
 and set $\bP= \underline{Isom}(\zeta\oplus \zeta^{-1},\cO^{\oplus 2}_X)\times_{\varrho}G$.
 Assuming all the Lie-algebraic data from Section \ref{Lie} fixed, we equip $\bP$
with the Higgs field $\theta= \dd \varrho(\theta_0)$, which can be identified with  the   matrix $y=\dd \varrho(y_0)=\sum_i f_i\in\fg$.
We shall discuss this identification in more detail in the next subsection.

 Specifying a  complex structure on $X$ is equivalent to specifying a conformal class of Riemannian metrics. A   metric $g$ within that  class
induces an hermitian metric on all tensor powers  $K_X^{\otimes m}$,  and more generally, on $\zeta^{\otimes m}= K^{m/2}_X$, $m\in\ZZ$, so we get
 a reduction of the structure group of $\bP$ to $U(1)=\varrho(U(1))\subset G$. If $\nabla$ is the corresponding
Chern connection, $F(\nabla)$ its curvature, and $F_1$ the curvature of the Levi-Civita connection, 
then Hitchin's equation
\[
 F(\nabla)+[\theta,\theta^\ast]=0
\]
reduces to ($\textrm{rk } \fg$ copies of) the equation $F_1-4i\omega_X=0$.
In other words, the  $U(1)$-reduction gives the harmonic metric for $(\bP,\theta)$
if and only if the Gauss curvature
$K_g=-4$. 
 This is shown for $G=SL(2,\CC)$ in \cite{hitchin_sd}. The  extension to other
groups is trivial and will be clear from the discussion that follows. It can also be deduced from the
functoriality (with respect to $G$) of non-abelian Hodge theory.
From the works of Poincar\'e and Koebe it is known that there is a unique
such metric in a given conformal class: it descends from the 	 standard hyperbolic metric on the upper half-plane
 after identifying the latter (biholomorphically) with the universal cover, $\widetilde{X}$, of $X$.
In this sense the harmonic (Hermite-Yang-Mills) metric on $(\bP,\theta)$ ``is'' the uniformising metric, and
 we call $(\bP,\theta)$ \emph{the uniformising Higgs bundle}, following Simpson
(\cite{simpson_uniformisation}).

The choice of a Killing form and a compact real form determine an hermitian product on $\fg$  (see Section \ref{Lie}). 
The harmonic reduction of  $\bP$   gives rise to an hermitian inner product on 
$\ad\bP$, which  is the harmonic metric for the Higgs (vector) bundle $(\ad\bP,\ad\theta)$.
 We also get  $L^2$-inner products on $A^p(\ad\bP\otimes K_X)$ for various $p$.

The infinitesimal deformations of the uniformising Higgs bundle (as well as those of any Higgs bundle, \cite{Biswas-Ramanan})  are  computed by the  Dolbeault complex
$$\xymatrix@1 { \ad \bP \ar^-{\ad\ \theta}[r]\ & \ \ad \bP\sotimes K_X}. $$
Taking its
Dolbeault resolution and passing to global sections we obtain the double complex

 $$\xymatrix{A^{0,1}(\ad \bP )\ar[r]^-{-\ad \theta} & A^{1,1}(\ad\bP)\\
A^{0,0}(\ad\bP)\ar[r]^-{\ad \theta}\ar[u]^{\dbar} & A^{1,0}(\ad\bP )\ar[u]^{\dbar}   },$$
whose total complex   is  
 \begin{equation}\label{dolbeault}
 \xymatrix@1{ 0\ar[r]&\ A^0(\ad\bP)\ar^-{d_0}[r]&\ A^{1,0}(\ad \bP)\oplus A^{0,1}(\ad\bP)\ar^-{d_1}[r]&\ A^{1,1}(\ad \bP)\ar[r]&\ 0}
 \end{equation}
with differentials $d_0=\left( \begin{array}{c}
                                        \ad \theta\\
                                        \dbar_{\bP}\\
                                \end{array}\right)$, $d_1= \left( \begin{array}{cc}
                                                                        \dbar_{\bP}, &- \ad\ \theta\\
                                                                   \end{array}\right)$. 

The dgla controlling the deformations of
$[(\bP,\theta)]\in M_{Dol}(G)$ is the  deformation complex (\ref{dolbeault}), i.e.,
$L^\bullet=A^\bullet(\ad \bP)$, with $d= \dbar_\bP + \ad \theta$ and 
 the standard bracket.
Notice that one can think of $\theta$ either as a twisted  section of $\ad\bP$, or as 1-form with values
in $\ad\bP$, and alternating between the two viewpoints may cause sign changes.
 The complex (\ref{dolbeault}) is
a slightly generalised version of Example \ref{higgs_dgla}, and it is immediate to see that conditions (\ref{condition1}) and (\ref{condition2}) from
Section \ref{symplectic_kur} are satisfied.
The Maurer-Cartan equation is
\begin{equation}
\label{MC1}
\dbar_{\bP} h + [ \theta + h,v ]  = 0, 
\end{equation}
$(h,v)\in A^{1,0}(\ad \bP)\oplus A^{0,1}(\ad\bP)$.
One  sees immediately  that if $(h,v)\in\MC(L^\bullet)$ and  $h$ is  holomorphic for $\dbar_{\bP}$, then
 $v\in\fz(\theta+h)$. This suggests that we can use the results and setup from
Sections \ref{symplectic_kur} and  \ref{centralisers}.
For that, we have to identify harmonic representatives of $H^1(L^\bullet)$ and see if condition (\hyperref[condition3]{3})   from Section \ref{symplectic_kur} is satisfied.
First  we discuss the structure of $\ad\bP$ in more detail.

\subsection{Filtrations, gradings and adjoints}
The homomorphisms between filtered (graded) objects in an abelian category carry a filtration (grading), and hence
the principal gradings on $\fg$ and  $\ad\bP$  induce gradings on their respective endomorphisms.
In particular, we have $\ad\in \mhom_0(\fg,\mend({\fg}))$, i.e., $\ad\in\mhom(\fg_m,\textrm{End}_m(\fg))$ for all $m$.
For the adjoint bundle and its endomorphism bundle we have
\[
 \ad \bP = \bigoplus_{m=-\cih}^\cih \ad_m\bP = \bigoplus_{m=-\cih}^\cih  \fg_m\otimes_\CC K_X^{\otimes m},
\]

\[
\send(\ad\bP)= \bigoplus_{m=-\cih}^\cih \send_m(\ad\bP) = \bigoplus_{m=-\cih}^\cih   \textrm{End}_m(\fg)\ctimes  K^{\otimes m}_X,
\]
and \emph{mut.mut.} for $\cA^{p,q}(\ad\bP)$ and $\cA^{p,q}(\send (\ad\bP))$.

Tensoring  with  powers of $K_X$ we  obtain plenty of trivial bundles: for all $m\in\ZZ$, 
\[\ad_m\bP\sotimes K_X^{-m}=\fg_m\ctimes\cO_X,
\]
\[
\send_m(\ad\bP)\sotimes K_X^{-m}= \textrm{End}_m(\fg)\ctimes \cO_X. \]
For every $m\in\ZZ$, $K^m\otimes_{\cO_X} K^{-m}$ has a canonical section $1_m$, namely, the image of $1\in \CC= H^0(X,\cO_X)$ under
$\cO_X\simeq K^m\otimes_{\cO_X} K^{-m}$ and  we have a commutative diagramme

\[
\scalebox{0.95}{
\xymatrix{
 \fg_{m}\ar[r]^-{\otimes 1_m}\ar[d]^-{\textrm{ad}}&\Gamma(X,\ad_m\bP\otimes K_X^{-m})\ar[d]^-{\textrm{ad}}\ar[dr]& \\
\textrm{End}_m\fg\ar[r]^-{\otimes 1_m}&\Gamma(X, \send_m(\ad\bP)\otimes K_X^{-m})\ar[r]^-{\io}&\mhom_{C^\infty} (A^\bullet(\ad\bP),A^\bullet(\ad\bP\otimes K_X^{-m})) . }}
\]
In particular,  for every $m\in \ZZ$ there are   inclusions
\[
\fg_m\hookr \mend_m \fg\hookr \Gamma(\send_m(\ad\bP)\otimes K_X^{-m})\hookr \mhom_{C^\infty} (A^{\bullet}(\ad\bP),A^{\bullet}(\ad\bP\otimes K_X^{-m})),
\]
\[
 \fg_m\ni \lambda\mapsto \io (\ad_\lambda\ctimes 1_m)\in \mhom_{C^\infty} (A^{\bullet}(\ad\bP),A^{\bullet}(\ad\bP\otimes K_X^{-m})) .
\]
For readability, we may occasionally suppress $\io$ or the subscript $m$ in $1_m$.

We get similar inclusions  if we   fix a   K\"ahler metric $h\in A^{0,1}(K_X)$ with  K\"ahler form $\omega_X\in A_X^{1,1}$:
\[
 \xymatrix@1{ \fg_m\ar[r]^-{\otimes \omega_m} & A_X^{1,1}(\ad_m\bP\otimes K_X^{-m}) }
\]
and
\[
 \xymatrix@1{ \fg_m\ar[r]^-{\otimes h _m} & A_X^{0,1}(\ad_m\bP\otimes K_X^{-m+1}) }.
\]

The natural  isomorphism $\cA^0 = \cA^{1,0}(K^{-1}_X)$ gives rise to a ``shift isomorphism'' 
$s: A^{\bullet,\bullet}(\ad\bP\sotimes K_X^{-m})\simeq A^{\bullet+1,\bullet}(\ad\bP\sotimes K_X^{-m-1}) $.
Again, we are going to suppress $s$ occasionally, but	 one should keep
in mind that
for $S\in \mend_{-1}\fg$, $\dbar (s\io S) + (s\io S)\dbar =0$, in particular, $\dbar$ 
\emph{anti-commutes} with $\ad_\lambda$. This is a consequence of \cite{voisin1},
 Remark 5.11 : the $\dbar$ operators on $\cA^{p,q}$ and $\cA^{0,q}(K^{\otimes p})$ differ by  
 $(-1)^p$.

We  now make some comments on 
adjoints and Hodge stars in order to  clarify conventions.

We are going to denote the hermitian metric on $T_X$ by $h$. In a local chart $(U,z)$ it is given  by
$h=\bh dz\otimes \dzbar$, and the K\"ahler form is  $\omega_X=\frac{i}{2}\bh dz\wedge \dzbar$.
The Riemannian metric $g$ on $T_{X,\RR}$ can be extended \emph{sesqui-linearly} to an hermitian pairing $g_{\overline{\CC}}$ on
$T_{X,\CC}$, which can then be restricted to $T_X^{1,0}$. Similarly for $T_{X,\CC}^\vee$ and its exterior powers.
 The pairing on $T_{X,\CC}^\vee$  equals \emph{half} of the direct sum of hermitian metrics on $\cA^{1,0}\oplus \cA^{0,1}$:
 see for example, \cite{voisin1},  Lemma 5.6
or \cite{huybrechts_cg}, Lemma 1.2.17. The Riemannian metric $g$ induces a dual metric, $g^\vee$ on $T_{X,\RR}^\vee$,
and, consequently, hermitian metrics $\widetilde{h}$ (on $K_X$) and $(g^\vee)_{\overline{\CC}}$ (on $T^\vee_{X,\CC}$).
One can check easily that   $(g^\vee)_{\overline{\CC}}=(g_{\overline{\CC}})^\vee$. However,
$\widetilde{h}= 4h^\vee$, where $h^\vee=\bh^{-1}\partial_z\otimes\partial_{\zbar}$ is the dual metric to $h$.

We are going to use the convention that the  Hodge star is \emph{anti-linear}, $\ast: A^{p,q}\to A^{1-p,1-q}$,
satisfying 
$\beta\wedge\ast\alpha= g(\beta,\alpha)_{\overline{\CC}\ } \omega_X$.  
On $1$-forms   $\ast$ coincides with conjugation up to  $\pm i$:
we have $\ast\alpha = i\overline{\alpha}$ for $\alpha\in A^{1,0}$.
     An  hermitian bundle, $E$, comes with an  anti-linear isomorphism
$\# :E\simeq E^\vee$, $e\mapsto \langle\ ,e \rangle$, where $\langle\ ,\rangle$ is the hermitian metric.
Notice that for $\alpha\in \Gamma_U (T^\vee_{X,\CC})$, $\#\alpha = -\frac{\ast \alpha\wedge}{\omega_X}$.
 We extend  $\ast$ by $\#$ and define  
$\ast: A^{p,q}(E)\to A^{1-p,1-q}(E^\vee)$ by
$\ast(\alpha\otimes e)=  \ast\alpha\otimes e^{\#}$. 

Let  $L$ and $M$ be hermitian line bundles on $X$,  $U\subset X$ a trivialising analytic open set,
 and $\lambda\in H^{0}(U,L)\simeq \mhom_U(M, M\otimes L)$ a nowhere vanishing section.
Then it is immediate to check that 
\[
 \lambda^\ast = \|\lambda\|^2 \lambda^\vee =\#\lambda.
\]
Here $\lambda^\ast\in A^0(\shom(L\otimes M, M))$ is the hermitian adjoint of $\lambda$ and
 $\lambda^\vee=\lambda^{-1}\in \mhom_U(L\otimes M,M)$ is the unique section pairing to $1$ with $\lambda$.
In particular, for a (nonvanishing) section $\lambda\in H^{0}(U,K_X)$ we have $\#\lambda=\frac{1}{2}\lambda^\ast$

 Let  $^\ast$ be a real structure on a vector bundle $E$ (compatible with the hermitian structure).
We extend $^\ast$ to $A^{p,q}(E)$ by \emph{complex conjugation}:
$(\alpha\otimes v)^\ast= \overline{\alpha}\otimes v^\ast$.
In particular, if $(U,z)$ is a local chart on $X$ and
  $\theta=\theta_z dz\in A^{1,0}(U,\ad\bP)$, 
 we have  $\theta^\ast= \theta_z^\ast d\zbar = i\theta_z^\ast \ast dz$.
This agrees with the conventions in \cite{hitchin_sd}; in \cite{hbls} the same quantity is denoted by $\overline{\theta}$.

As a special case, let us consider $1\in H^0(\cO_X)\subset A^{1,0}(\shom(M,M\otimes K^{-1}_X))$,
where $M$ is an arbitrary hermitian line bundle. It is immediate to check that
$1^\ast = h = \bh dz\otimes \dzbar\in A^{0,1}(\shom(M\otimes K^{-1}_X,M))$.
Here $1^\ast$ means, naturally, $(s(1))^\ast= (1_1)^\ast$.
More generally, given $\lambda\in \fg_{-1}$, $\lambda\otimes 1 = s(\lambda\ctimes 1_{-1})\in A^{1,0}(X,\ad_{-1}\bP)$ and we have  
$(\lambda\otimes 1)^\ast =  \lambda^\ast\otimes h\in A^{0,1}(\ad_1\bP)$. Similarly
$  (\ad_\lambda \otimes 1 )^\ast= (\ad_{\lambda^*}\otimes h)\in A^{0,1}(\send_1\ad\bP)$.  If we have two elements
$\mu,\lambda\in \fg_{-1}$, then $[\mu\otimes 1,(\lambda\otimes 1)^\ast]=-2i[\mu,\lambda^\ast]\ctimes \omega_X$.

Finally, we can consider $\ad_\lambda\ctimes 1$ as an operator acting on $A^\bullet(\ad\bP)$, in which case its  adjoint then is
$\frac{1}{2}\ad_{\lambda^\ast}\ctimes 1$. More pedantically, 
$\left( \io s\ \ad_\lambda\ctimes 1_{-1} \right)^\ast = \frac{1}{2}\io \ad_{\lambda^\ast}\ctimes 1_{1}$.

\subsection{Harmonic Representatives of cohomology}\label{harmonic_reps_coho}
Now we return to the Dolbeault complex (\ref{dolbeault}).  We have

 \begin{equation}\label{dolbeault_triv}
\scalebox{0.85}{
 \xymatrix@1{ 0\ar[r]&\ \bigoplus_mA^0(\ad_m\bP)\ar^-{d_0}[r]&\ \bigoplus_m A^{1,0}(\ad_m\bP)\oplus  A^{0,1}(\ad_m\bP)\ar^-{d_1}[r]&\ \bigoplus_mA^{1,1}(\ad_m\bP)\ar[r]&\ 0},}
 \end{equation}
with differentials 
 $d_0=\left( \begin{array}{c}
                                         \ad_y\otimes 1\\
                                        \dbar\\
                                \end{array}\right)$ and $d_1= \left( \begin{array}{cc}
                                                                        \dbar, & \ad_y\otimes 1\\
                                                                   \end{array}\right)$.
For legibility, we have suppressed the obvious part of  the nomenclature: $\ad_y\otimes 1$ stands for $\io \ad\theta= \io s\ \ad_y\otimes 1_{-1}$
and 
the  Dolbeault operator
$\dbar$ is $\oplus_m \dbar_{K^m}$,  the direct sum of the Dolbeault operators on $\ad_m\bP=\fg_m\ctimes K^{\otimes m}_X$.

\begin{thm}
\label{H1}
Let $\fg=\textrm{Lie}(G)$ be a simple complex Lie algebra, equipped with $\{x,2\rho^\vee,y\}$ as
in Section \ref{Lie}, so $\fz(x)\simeq \oplus_i \fg_{m_i,i}$ and $\fz(y)=\oplus_i \fg_{-m_i,i}$, where
$\fg_{\pm m_i,i}=\fg_{\pm m_i}\cap W_i$.
Let $(\bP,\theta)$ be the 
uniformising $G$-Higgs bundle, equipped with the Hermite-Yang-Mills metric and
let $L^\bullet$ be the Dolbeault complex
(\ref{dolbeault_triv}). Then
\[
 H^1(L^\bullet)\simeq \cH^1(L^\bullet) \subset A^{1,0}(\ad\bP)\oplus A^{0,1}(\ad\bP),
\]
where
\[ 
 \cH^1(L^\bullet)=  \oplus_i\fg_{m_i,i}\ctimes\cH^{1,0}( K_X^{m_i})\bigoplus \oplus_i\fg_{-m_i,i}\ctimes\cH^{0,1}( K^{-m_i}_X).   
\]
Hence 
\[
 \cH^1(L^\bullet)\simeq  \cH^{1,0}(\fz(x)_\bP)\oplus \cH^{0,1}(\fz(y)_\bP)\simeq \cB_\fg\oplus \cB_\fg^\vee,
\]
where   $\cB_\fg=H^0(X,\ft \otimes K_X/W)$ denotes  the Hitchin base.
\end{thm}
\begin{remark}
The vector bundles $\fz(x)_\bP$ and $\fz(y)_\bP$ are  (the obvious) twists of the centralisers $\fz(x)$ and $\fz(y)$ by $\bP$:
$\fz(x)_\bP = \oplus \fg_{m_i,i}\ctimes K_X^{m_i}$ and $\fz(y)_\bP = \oplus \fg_{-m_i,i}\ctimes K_X^{-m_i}$. Since  $K_X^{m}$
are hermitian, we can talk about harmonic represenatives of their cohomology, hence the notation. Explicitly, 
$\cH^{1,0}(\fz(x)_\bP)= \oplus_i \fg_{m_i,i}\ctimes \cH^{1,0}(X,K_X^{m_i})$ and $\cH^{0,1}(\fz(y)_\bP)= \oplus_i \fg_{-m_i,i}\ctimes \cH^{0,1}(X,K_X^{-m_i})$.
 Recall  that   $\dim \fg_{\pm m_i,i}=1$. We are assuming a fixed basis for the $G$-invariant polynomials on $\fg$, which fixes, by duality, bases for 
the spaces $\fg_{m_i,i}$, as discussed in Section \ref{Lie}. 
These give rise to bases of  $\fg_{-m_i,i}$: either  by taking  hermitian conjugates  or by applying $\ad_y^{2m_i}$ (the two choices differ by a combinatorial coefficient). 
The identification $\cB_\fg \simeq  H^0(X,\bigoplus_i K_X^{m_i+1})$ depends on the choice of invariants polynomials.
The identification  $\cH^1(L^\bullet)\simeq  \cH^{1,0}(\fz(x)_\bP)\oplus \cH^{0,1}(\fz(y)_\bP)$
depends on the choice of basis for $\fg_{\pm m_i,i}$ and uses the hermitian metric.
\end{remark}

\emph{Proof:}\\
For the purposes of the proof, let us denote the summands in  the decomposition 
\[\fg =\fz(x)\oplus \left( \textrm{Im}(\ad_x)\cap \textrm{Im}(\ad_y)\right) \oplus\fz(y)
\]
by subscripts $x$, $o$ and $y$, so $\fg=\fg_{x}\oplus \fg_{o}\oplus \fg_y$. Use  combinations of subscripts   to denote projections on
pairs of  summands.
If $\sigma=(\sigma',\sigma'')^T\in \ker d_1\subset A^1(\ad\bP)$, then  
\[
 \sigma= d_0(P\otimes 1 (\sigma '_{oy})) + (\sigma '_x,0)^T + (0,\sigma_{y}'')^T, \ \dbar \sigma '_x=0.
\]
The first summand  is a coboundary and the second term is never a coboundary, as $\fz(x)\simeq \cok(ad_y)$. The  last summand, however, can contain a 
$\dbar$-exact term. By the  Hodge decomposition on $A^p(X,K^m_X)$, we can write $\sigma''_{y}\in A^{0,1}  (\fz(y)_\bP)$ as
$\sigma''_{y}= \dbar\left(\dbar^*\bG \sigma''_{y}\right) + \bH(\sigma''_{y})$, where $\bG$ is Green's operator. Thus altogether
\[\sigma= (\sigma',\sigma'')^T = d_0\left(P\ctimes 1 (\sigma'_{oy}) + \dbar^*\bG \sigma''_{y} \right) + (\sigma'_x,\bH \sigma_{y}'')^T,\]
$\dbar \sigma'_x=0$, and we obtain
\[
 \ker d_1 =\textrm{Im}d_0 \bigoplus \oplus_i\fg_{m_i,i}\ctimes \cH^{1,0}( K_X^{m_i-1})\bigoplus \oplus_i\fg_{-m_i,i}\ctimes \cH^{0,1}( K^{-m_i+1}_X).
\]
The second and third direct summands are hence isomorphic to $H^1(L^\bullet)$, and are identified (via  shifts)
with  $\cB_\fg\oplus\cB^\vee_\fg$.
In the next proposition we show that these are actually the harmonic representatives for $H^1(L^\bullet)$.
 Explicitly, the isomorphism 
  $H^1(L^\bullet) \simeq  \cH^{1,0}(\oplus_i K_X^{m_i})\oplus \cH^{0,1}(\oplus_i K^{-m_i}_X)$
is given by
\[
 [\sigma]=[(\sigma',\sigma'')]\mapsto (\sigma'_x,\bH \sigma_{y}'').
\]

\qed
\begin{remark}
Using the explicit knowledge of the differentials of $L^\bullet$, one can check easily that $H^2(L^\bullet)=0=H^0(L^\bullet)$.
 Moreover,  $\textrm{Aut}(\bP,\theta)=Z(G)$, i.e., the pair  has no ``extra'' automorphisms (it is regularly stable). This can be deduced for instance from
Proposition 3.1.5 (ii), \cite{bei-drin} and the non-abelian Hodge theorem (\cite{hbls}). Hence $[(\bP,\theta)]$ corresponds to a smooth point of $M_{Dol}(G)$. 
Of course, this is already contained  in \cite{hitchin_teich} for the case when $G$ is of adjoint type.
\end{remark}

\begin{Prop}\label{harmonic}
 The vector space  $\cH^1(L^\bullet)=\cH^{1,0}(L^\bullet)\oplus \cH^{0,1}(L^\bullet)$ is the space of harmonic
representatives of $H^1(L^\bullet)$. That is, 
\[
 \cH^1(L^\bullet)= \ker d_1\cap \ker d_0^\ast\simeq H^1(L^\bullet).
\]
\end{Prop}

\emph{Proof:}\\
We have
\[d_0^* =\left(\ad^\ast_{\theta},\dbar^* \right) = \left(2\  s^{-1} \ad_x\ctimes 1,\dbar^* \right):A^{1,0}(\ad\bP)\oplus A^{0,1}(\ad\bP)\to A^0(\ad\bP),
\]
and   $\sigma\in \ker d_0^*$  implies
$$ \sigma=(\sigma',\sigma'')^T = (\sigma'_x,0)^T+ (0,\sigma''_{y})^T + \left(-\frac{1}{2}(Q\ctimes 1)(\dbar^*\sigma''_{xo}),\sigma''_{xo} \right)^T,\  \dbar^*\sigma''_{y}=0. $$
Applying the Hodge decompositions $A^{0,1}(K^{-m_i+1}_X)=\textrm{Im}\dbar\oplus\cH^{0,1}(K^{-m_i+1})$ and
$A^{1,0}(K^{m_i+1}_X)= \cH^{1,0}\oplus \textrm{Im}\dbar^\ast$ to $\sigma_y''$ and $\sigma'_x$, respectively, 
 we get
\[
 \ker d_0^\ast =\textrm{Im}d_1^\ast \bigoplus \oplus_i\fg_{m_i,i}\ctimes\cH^{1,0}( K_X^{m_i})\bigoplus \oplus_i\fg_{-m_i,i}\ctimes\cH^{0,1}( K^{-m_i}_X),
\]
and the result follows. 
\qed
\begin{remark}
For the  proofs of Theorem \ref{H1} and  Proposition \ref{harmonic} it is not essential that  the principal $\slt$  is  related 
in a specific way to some fixed Chevalley generators, but it is  essential that principal subalgebra is real,   the
compact anti-involution  maps $\fg_k$ to $\fg_{-k}$, and  the different $W_i$'s are mutually orthogonal. 
\end{remark}

\begin{Prop}
The induced complex symplectic form on 
\[
  \cH^1(L^\bullet) \subset A^{1,0}(\ad\bP)\oplus A^{0,1}(\ad\bP)
\]
is the
canonical symplectic form on
$\cH^{1,0}(\fz(x)_\bP)\oplus \cH^{0,1}(\fz(y)_\bP)$
and agrees, up to Lie-theoretic normalisation factors, with the canonical symplectic form on $\cB_\fg\oplus \cB^\vee_\fg$.
\end{Prop}

\emph{Proof:}\\
 This  is essentially clear from the construction. The Killing form $\kappa$ places $\fz(x)$ and $\fz(y)$ in duality. Next,
the complex symplectic form on $H^1(L^\bullet)$ is induced by 
 the (weak) duality pairing
\[\left( A^{1,0}(\ad\bP)\oplus A^{0,1}(\ad\bP)\right)^{\times 2}\to A^{1,1}\to \CC,\]
\[
 ((u,\alpha),(v,\beta))\mapsto \int_X \kappa(u\wedge\beta)) - \kappa(v\wedge\alpha).
\]
Evaluating it on  pairs of harmonic representatives of $H^0(K^{m_i+1})$ and $H^1(T^{m_i}_X)$, say,
$ (u_i,\alpha_i)$, $(v_i, \beta_i)$, we get an expression of the form
$\sum_i \kappa(e_{m_i},e_{-m_i})(\beta_i(u_i)- \alpha_i(v_i))$, where $e_{m_i}$ and $e_{-m_i}$
are bases of the 1-dimensional vector spaces $\fg_{\pm m_i,i}$. If they are  dual bases, 
then the pairing will coincide with  the canoncial symplectic form on $\cB_\fg\oplus \cB_\fg^\vee$, otherwise
there will be extra coefficients $\kappa(e_{m_i},e_{-m_i})$. \qed

\subsection{The Symplectic Kuranishi slice}
Now we  apply the results from the previous sections to the deformation theory of the uniformising Higgs bundle.
\begin{thm}\label{exponential}
 Keep the notation and assumptions from the previous sections. Consider the holomorphic family of Higgs bundles
\[
\Gamma: \cB_\fg\times\cB_\fg^\vee \longrightarrow  A^{1,0}(\ad\bP)\oplus A^{0,1}(\ad\bP),
\]
\[
 \Gamma(h,v)= \left(h,\Phi_h^{-1}(v)\right) = \left(h, \sum_{k=0}^{\cih}(-1)^k\left((s^{-1}P\ctimes 1)\circ ad_{h} \right)^k(v) \right),
\]
where 
\[
 (h,v)\in \oplus_i\fg_{m_i,i}\ctimes\cH^{1,0}( K_X^{m_i})\bigoplus \oplus_i\fg_{-m_i,i}\ctimes\cH^{0,1}( K^{-m_i}_X)\simeq \cB_\fg\times\cB_\fg^\vee
\]
and
\[\Phi_h= 1+s^{-1}(P\ctimes 1)\ad_h\in \mend (A^{0,1}(\ad\bP)).\]
The family $\Gamma$
is a miniversal deformation of the uniformising Higgs bundle $(\bP,\theta)$. There exists an open neighbourhood  $\cU\subset\cB_\fg\times\cB_\fg^\vee$
containing $0$,  
for which $\Gamma\vert_{\cU}$ is a universal deformation. Moreover, $\Gamma^\ast\omega_{can}= \omega_{can}$.
\end{thm}
\begin{remark}
 Clearly, we also have a formal version of $\Gamma$, i.e.,  
 a functor of Artin rings $\underline{\Gamma}: \underline{\cB}_\fg\times\underline{\cB}^\vee_\fg\to \Def_L$, given by the same formula as above.
As everywhere above, $\cB_\fg$ should be understood in terms of harmonic representatives.
\end{remark}

\emph{Proof:}\\
From subsection \ref{harmonic_reps_coho} we know that the dgla $L^\bullet$ satisfies conditions
$(1)$,$(2)$ and $(3)$ from Section \ref{symplectic_kur}, and  that $H^2(L^\bullet)=0=L^3$.
In the notation of Section \ref{symplectic_kur}, $d_1'=ad\theta$, and we have a splitting,
$s^{-1}P\ctimes 1$.
By Lemma \ref{induction}, the geometric series for $\Phi_h^{-1}$ reduces to the given formula, i.e., $\pi$ drops out of the
expressions.
The condition $(\pi-1)[\Gamma,\Gamma]=0$ holds (essentially) for the same reasons as in Section \ref{centralisers}: since $h$ (resp. $v$)
is a highest (resp. lowest) weight vector, none of the sections from $[\Gamma,\Gamma]$ will be contained in $A^{1,1}(\fz(y)_{\bP})$.
The Kodaira-Spencer map of this family is the identity, so, by \cite{fukaya_def}, Theorem 1.3.3. this is a miniversal family.
Since $(\bP,\theta)$ is regularly stable, we obtain a universal family by restricting the domain of $\Gamma$.

Finally, the statement that $\omega_{can}$ pulls back to $\omega_{can}$ follows from Theorem \ref{triv_symp}. The conditions in that theorem are
satisfied (essentially) for the same reason as before: the pairing $A^{1,0}(\ad\bP)\times A^{0,1}(\ad\bP)\to \CC$ is obtained by combining
cup product  $A^{1,0}(K^{m_i})\times A^{0,1}(K^{-m_i})\to A^{1,1}\to \CC$  with the Killing form $\kappa:\fg\times\fg\to\CC$.
 But the Killing form is non-zero only on
$\fg_{m_i,i}\times \fg_{-m_i,i}$ and $\fg_{-m_i,i}\times \fg_{m_i,i}$, and  since  $s^{-1}(P\otimes 1) \pi\ad_h$ has degree 2 (with respect to the principal grading), the 
orthogonality condition 
 from \ref{triv_symp} is satisfied.
 \qed 

\begin{remark}
 The description of the family from  Theorem \ref{exponential} is constructive, and one can write $\Gamma$ explicitly once the Lie-algebraic data are fixed. After all,
$\Gamma$ is simply a ``twisted'' version of the analogous formula from Section \ref{centralisers} (see the example there.)
\end{remark}
We now rephrase the result and draw some easy corollaries.
\begin{Cor}
Denote, as in Section \ref{symplectic_kur}, $\cS=\MC(L)\cap\ker[(1-\bH')\pi']$. Then
 $\Phi=\Gamma^{-1}:\cS\to \cB_\fg\times \cB_\fg^\vee$ provides Darboux coordinates on $\cS$.
\end{Cor}
Notice, once again, that if we want to consider $\cS$ with its (somewhat useless) structure of a (germ of a) subvariety
of an infinite-dimensional vector space, we should first complete $L^\bullet$ with respect to a suitable
(H\"older or Sobolev) norm. However, $\Gamma(\cB_\fg\times\cB^\vee_\fg)\subset A^1(\ad\bP)$ since 
$s^{-1}(P\ctimes 1)\pi\ad_h$ is nilpotent. See also Corollary \ref{nilpotent_cor}.

In \cite{hitchin_teich}, N.Hitchin constructed a section, $\fs:\cB_\fg\to M_{Dol}(G)$ as  a ``global version'' of 
Kostant's section $\ft/W\to \fg$.
We have identified $\cB_\fg\simeq\oplus_i \cH^{1,0}(X,K^{m_i}_X)$, and have embedded the latter
into $A^{1,0}(\ad\bP)$ via the basis vectors $e_{m_i}$, spanning $\fz(x)$. The section then is 
the holomorphic family of Higgs bundles, whose
underlying bundle is $\bP$, and which carries  the Higgs field $\theta+\sum_i e_{m_i} \alpha_i$,
$\alpha_i\in\cH^{1,0}(X,K^{m_i}_X)$ .
In terms of deformation functors, the section is given by $\fs:\underline{\cB}_\fg\to A^{1}(\ad\bP)\to\Def_L$, 
$\fs(h)=(h,0)$.
\begin{remark}
 What we give here is  a somewhat non-canonical description of the section.  To construct the section, 
one only needs a choice of theta-characteristic, $\zeta$, and a principal $\slt$-subalgebra.
 For a more detailed and geometric treatment,
see \cite{hitchin_teich}, \cite{don-pan}, \cite{ngo}.
\end{remark}

\begin{Cor}
The restriction of 
 $\Gamma$ to $\cB_\fg\times\{0\}$ is the  Hitchin section. If we  regard elements of
$\cB^\vee_\fg\simeq \oplus \cH^{0,1}(K^{-m_i}_X)$ as linear Hamiltonian functions on the base, we get  that
$\Gamma(h,v)=\exp_{X_v}(\fs(h))$,  where $X_v$ is the Hamiltonian vector field, corresponding to $v$.
\end{Cor}
\emph{Proof:}
 We have $\Gamma(h,0)= (h,0)=\fs(h)$ by construction. The rest is immediate from the Theorem. \qed

 Because of the above, one may  refer to $\Gamma$ as a ``holomorphic exponential map''. 

\begin{remark}
We can also look at the image of $\{0\}\times\cB^\vee_\fg$ under $\Gamma$. 
This is a family of Higgs bundles for which the Higgs field is constant (as a smooth twisted endomorphism), while the holomorphic structure
varies in $\cB_\fg^\vee$. For instance, if 
 $\fg=\slt$, the underlying vector bundles are  extensions of $\zeta$ by $\zeta^{-1}$. 
They all come with a canonical inclusion $\CC\hookr H^0(\send E\otimes K_X)$, and the Higgs field is the image of $1\in\CC$.
\end{remark}

Recall from Section \ref{centralisers} that $ I\sslash G\simeq I\vert_{y+\fz(x)}= \bs^*I$, where $I$ is the tautological family of centralisers
and coincides with the set of Maurer-Cartan elements for a certain dgla. There is a ``global version'' of this statement.
More precisely, 
recall that a  ``regular Higgs field'' is one which is an everywhere regular section of $\ad\bP\otimes K$, i.e., pointwise it takes values in
 $\fg^{reg}$.  
\begin{Cor}
The image of the exponential map consists of all regular Higgs bundles in the connected component
of the uniformising Higgs bundle.
\end{Cor}
\emph{Proof:}\\
Since the values of  $\Gamma$ are regular by construction, 
 the nontrivial statement is the opposite  inclusion. We claim that
every Higgs pair with  a regular Higgs field is isomorphic to one in the image of $\Gamma$.
 Suppose $(\bQ,\fii)$ is such a Higgs bundle.
Then there exists a  $C^\infty$-isomorphism 
$\ad\bQ\sotimes K_X\simeq_{C^\infty}\ad\bP\sotimes K_X$, and
$\dbar_\bQ= \dbar_\bP +\xi$. Since $\fii$ is regular, it can be  conjugated
to $\fs(\chi(\fii))$. This replaces  $\dbar_\bQ$  by a gauge-equivalent Dolbeault operator, 
$\dbar_\bP+ v$, 
and the Maurer-Cartan equation states that $v\in\fz(\fs(\chi(\fii)))$, i.e., $(\fs(\chi(\fii)),v)\in\cS\simeq \cH^1(L^\bullet)$.
Notice that we are not making any claim regarding stability of our bundles.
\qed
\begin{remark}
 In unpublished notes (\cite{teleman_langlands}) C.Teleman proved the same result for $GL(n,\CC)$.
\end{remark}

\section{Appendix: Kuranishi Theory}\label{Kuranishi_Theory}

In this subsection we recall some relevant facts from formal and analytic Kuranishi theory. 
Our main references will be
\cite{gm_kur}, \cite{goldman-millson} and \cite{maurer}, and, to a lesser extent, 
\cite{fukaya_def} and \cite{kontsevich_def}.

Suppose $L^\bullet$ is a dgla equipped with a \emph{splitting}
$\delta$. This is a linear map $\delta \in \mhom^{-1}(L^\bullet,L^\bullet)$  which satisfies $\delta^2=0$, $d=d\delta d$ and $\delta= \delta d \delta$.
Notice that while $d$ is a derivation of the bracket, $\delta$
\emph{a priori} need not be compatible with the bracket in any way. In fact,  if $L^\bullet$ admits a splitting 
which is a derivation, then it  is \emph{formal}
(\cite{kosarew}, Theorem 4.2.1.).  For comparison, any choice of splitting gives an isomorphism between $L^\bullet$ and 
$H^\bullet(L)$ as $L_\infty$-algebras: see e.g., \cite{kontsevich_def} or \cite{fukaya_def}.

Specifying a  splitting is equivalent to specifying a
 ``Hodge decomposition'' $L^\bullet = \cB^\bullet\oplus \cH^\bullet \oplus\cC^\bullet$, where
$\cB^\bullet = \img (d)$, $\cC^\bullet = \img(\delta)$, $\cH^\bullet= \ker d\cap \ker \delta\simeq H^\bullet(L^\bullet)$
and  $\cH^\bullet\oplus\cC^\bullet =\ker \delta$. 
The
map $\delta$  gives a (co)chain homotopy 
 between  $\bH=\textrm{pr}_{\cH}$ and the identity, i.e., 
$d\delta +\delta d =1-\bH$. This can also be written as 
$\textrm{pr}_{\cB} + \textrm{pr}_{\cH} + \textrm{pr}_{\cC} = 1$
since
$d\delta = pr_{\cB}$ and $\delta d= pr_{\cC}$.
The choice of  splitting gives a decomposition of   $(L^\bullet,d)$ into a sum of 1-term complexes
$\cH^i[-i]$ and 2-term contractible complexes  
$\xymatrix@1{\cC^i\ar[r]^-{d}&\cB^{i+1}}$. 

For a given $\psi\in \cB^\bullet$, the equation  $d\fii=\psi$ has unique solution, $\fii = \delta\psi$ \emph{in} $\cC[-1]^\bullet$.
However, if  $\cH^\bullet\neq (0)$, this equation will have infinitely many solutions
in   $\ker \delta[-1]=(\cH^\bullet\oplus\cC^\bullet)[-1]$, since  $\fii = \bH(\fii) +\delta\psi$, and the harmonic part $\bH(\fii)$ can be arbitrary.

One 
 approaches the  Maurer-Cartan equation $d\fii = - \frac{1}{2}[\fii,\fii]$ in a  similar way: we look for all
$\fii\in L^1$
such that  $\fii = \bH(\fii) - \frac{1}{2}\delta[\fii,\fii]$.  
 They constitute   the zero locus of 
\[
 M:L^1\to \cB^1\oplus\cC^1
\]
\[
 M(\fii) = (1-\bH)\left(\fii+\frac{1}{2}\delta[\fii,\fii]\right) = (1-\bH)(F(\fii)),
\]
where $F(\fii)= \fii+\frac{1}{2}\delta [\fii,\fii]$ is 
the \emph{Kuranishi map} $F:L^1\to L^1$.
So
\begin{equation} \label{slice}
\ker\delta =\cH^1\oplus \cC^1\supset \left\{\fii:\ \fii = \bH(\fii) - \frac{1}{2}\delta[\fii,\fii] \right\}=M^{-1}(0)= F^{-1}(\cH^1).
 \end{equation}
\emph{A prirori}  $\fii\in F^{-1}(\cH^1)$ is \emph{not} a Maurer-Cartan element.
There is, however,  an obvious necessary condition that Maurer-Cartan elements have to satisfy:
  $ d\fii = - \frac{1}{2}[\fii,\fii] \Rightarrow \bH[\fii,\fii]=0$. Hence we   define 
$k:L^1\to\cH^2$ by  $k(\fii)=\bH[\fii,\fii]$ and look at  the set
 $ F^{-1}(\cH^1)\cap k^{-1}(0)$ and at its image under $F$,    $K_L := F(F^{-1}(\cH^1)\cap k^{-1}(0))\subset \cH^1$.

Loosely speaking,   if we work formally (or analytically), then (the germ of)  $ F^{-1}(\cH^1)\cap k^{-1}(0)$ 
consists of all  Maurer-Cartan elements
in (some neighbourhood of zero in) $\cH^1\oplus\cC^1$ and provides a semi-universal family of deformations of $\Def_L(\CC)$.

We first make some remarks about the Kuranishi map.   Since $\delta(F(x))=\delta(x)$, we have that
$F(\ker\delta)\subset \ker\delta$, and, by (\ref{slice}), $F^{-1}(\cH^1)\subset \ker\delta$.
Next,  
the ``slice'' $Y_L:=\cQ_L^{-1}(0)\cap\ker\delta$, consisting of Maurer-Cartan elements in $\ker\delta$,
gets mapped to $\cH^1$ by $F$.
 Indeed, 
$y\in Y_L \Rightarrow dF(y)=\frac{1}{2}\bH [y,y]\in \cB^2\cap\cH^2=(0)$, so
$F(Y_L)\subset \cH^1$. Moreover,
 $Y_L= F_{\vert Y}^{-1}(\cH^1)\subset M^{-1}(0)\cap k^{-1}(0)$, so $F(Y_L)\subset K_L$.
Notice that $F$, considered as a quadratic map between vector spaces  
(or subsets thereof)
need \emph{not} be invertible!

The next diagramme illustrates the different inclusions:
\[
 \xymatrix{Y_L=\MC(L)\cap\ker\delta\ar@{^{(}->}[r]\ar[d]&F^{-1}(\cH^1)\cap k^{-1}(0)\ar@{^{(}->}[r]\ar[d]&F^{-1}(\cH^1)\ar@{^{(}->}[r]\ar[d]&\ker\delta\ar@{^{(}->}[r]\ar[d]&L^1\ar[d]^-{F}\\
	    F(Y_L)\ar@{^{(}->}[r]&K_L\ar@{^{(}->}[r]&\cH^1\ar@{^{(}->}[r]&\ker\delta\ar@{^{(}->}[r]&L^1\\
}.
\]

In order to say more, we need a topology.

First we turn to  the formal setup and  define  a functor
$\YY_L=  \MC_L\cap\ker\delta\in\fart$,
$$ \mathbb{Y}_L(A)=Y_{L\otimes\fm_A} =\left\{ \eta\in \ker \delta\otimes \fm_A: d\eta+\frac{1}{2}[\eta,\eta]=0 \right\}.$$

For every element in $\MC_L(A)$ there is a unique gauge transformation, taking it to
$\YY_L(A)$, see e.g. \cite{simpson_eyssidieux}, Lemma 2.6. This is a variant of the so-called
Uhlenbeck slice (Coulomb gauge). Slice theorems have been  widely used in gauge theory since the late 1970's, 
most notably by    Atiyah--Hitchin--Singer, Taubes and Uhlenbeck and before that by Parker and
Mitter--Viallet.

The Kuranishi map gives rise to  a  functor  $\FF\in \fun(\underline{L}^1,\underline{L}^1)$, given by the same
formula as before, \emph{mut. mut.}
One shows, using artinian induction, that $\FF$  is an \emph{isomorphism}, see e.g., \cite{gm_kur}, Lemma 3.1 or \cite{maurer}, Lemma 4.2.

The \emph{Kuranishi functor} $\KK_L\in\fart$ is defined as the kernel of $k\circ \FF^{-1}\in\fun(\underline{\cH}^1,\underline{\cH}^2)$:
$$\KK_L(A)= K_{L\otimes\fm_A}= \left\{ x : \bH([\FF_A^{-1}(x),\FF_A^{-1}(x)])=0 \right\} \subset \cH^1\otimes \fm_A . $$

Applying the earlier considerations to  $L\otimes\fm_A$, we see that 
$\FF\in\fun(\YY_L,\KK_L)$, and it is in fact  an isomorphism
(\cite{gm_kur}, Section 3 or \cite{maurer}, Proposition 4.6).
Then
\[\xymatrix@1{\KK_L\ar[r]^-{\FF^{-1}}& \YY_L \ar[r]& \Def_L} 
\]
is shown to be \'etale:    see \cite{gm_kur}, Section 3 or \cite{maurer},
Theorem 4.7 for further details. The functor $\YY_L$ is called a \emph{formal miniversal deformation} or
\emph{formal Kuranishi space}. 
The isomorphism class of $\YY_L$ is independent of the choice of $\delta$ and 
quasi-isomorphic  dgla's have isomorphic $\YY_L$'s (\cite{gm_kur}) It is clear  that $\YY_L$ is 
pro-representable, i.e., $\YY_L= h_R$, for a complete local algebra $R$. An explicit description of $R$ can
be obtained by fixing a basis of $\cH^1$, say $\{\eta_1,\ldots,\eta_d\}$, with dual basis $\{t_1,\ldots,t_d\}$,
and then taking $R=\CC\{t_1,\ldots,t_d\}/\cI$. The ideal  $\cI$ is generated  by the components of
$\bH[\sum_i t_i\eta_i,\sum_j t_j\eta_j]=0$ with respect to some basis of $\cH^2$.
If $L$ is formal, then $\KK_L$ is the quadratic cone in $\underline{\cH}^1$ determined by cup product.
 If $H^2(L)=0$, $\KK_L=\underline{\cH}^1$
and $R=\CC\{t_1,\ldots,t_d\}$. If $H^0(L)=0$, then $\Def_L$ is pro-representable.

If $L^\bullet$ itself  carries a topology  we can go beyond the formal level and exhibit a germ of a complex space
whose local ring completed at the origin prorepresents $\KK_L$.
Suppose $L^\bullet$ is an \emph{analytic dgla}
in the sense of \cite{gm_kur}.
This means that $L^\bullet$ is a \emph{normed dgla}
(i.e., for all $i\in\NN$ there is a norm $\| \|_i$ on $L^i$ with respect to which $d$ and $[ , ]$ are continuous) 
 and the completion
$\widehat{L}^\bullet$ is equipped with a continuous splitting, $\delta$.
In the case of ``usual'' Hodge theory, the norms are the Sobolev norms and $\delta = \dbar^\ast \bG$,
where $\bG$ is Green's operator. In the case of deformations of a complex manifold or deformations of a holomorphic vector bundle the norms are H\"older norms.
The splitting $\delta$ has to be  compatible with the inclusion $L^\bullet\subset \widehat{L}^\bullet$, which  means two things.
 First, we assume that $\cH= \ker d\cap \ker\delta\subset L\subset \widehat{L}$.
And second,  we assume that the three projections preserve $L^\bullet\subset\widehat{L}^\bullet$
and $\textrm{pr}_{\cB }L=d(L[-1])$.
 Such a   $\delta$ induces
a splitting of $L^\bullet$ as well.

Then, by the implicit function theorem for Banach spaces
$F: \widehat{L}^1\to \widehat{L}^1$ is an analytic isomorphism between open balls around the origin: we have 
  $\dd F_\xi = 1 + \delta\ad_\xi$, and  $\dd F_0=1$ (see \cite{gm_kur}, Lemma 2.2).
We can introduce now the analytic versions of all of the above functors. First  set 
\[
\cY = Y_{\widehat{L}} =\left\{ \eta: \delta \eta=0, d\eta+\frac{1}{2}[\eta,\eta]=0 \right\} \subset \cH^1\oplus \widehat{\cC}^1.
\]
Notice that $\cY$  is an  algebraic subset of the   (possibly)
infinite-dimensional vector space $\ker\delta =\cH^1\oplus \widehat{\cC}^1$.  
Next, let
$$ \cK = K_{\widehat{L}} = \left\{ x \in \cH^1 : \bH([F^{-1}(x),F^{-1}(x)])=0\right\}\subset \cH^1.$$
The previous discussion, applied to $\widehat{L}$,  gives  $F(\cY)\subset \cK$.

The generalised  Kuranishi's theorem
(\cite{gm_kur}, Theorem 2.3) 
states that  the Kuranishi map  
induces an analytic isomorphism of germs $F: (\cY,0)\simeq (\cK,0)$, and
hence the functors $\KK_L$ and $\YY_L$ are prorepresented by $\widehat{\cO}_{(\cK,0)}$.
To prove it  one must show that for some open ball $B_0 \subset \widehat{L}^1$ we have
$F^{-1}(\cK\cap B_0) = B_0'\cap F^{-1}(\cH^1)\cap k^{-1}(0)\subset \cY$. 
Indeed, if $F(\xi)\in \cH^1$ and $\bH[\xi,\xi]=0$,  we get immediately that
$\delta F(\xi) =\delta(\xi)  0$.	
 Then
\[
 dF(\xi)=  d\xi + \frac{1}{2}d\delta[\xi,\xi] = 0 = \left( d\xi+\frac{1}{2}[\xi,\xi]\right) -\frac{1}{2}\delta d[\xi,\xi].
\]
We have to show that the last summand is zero.
The fact that $d$ is a derivation, combined with the Jacobi identity shows that $\delta d[\xi,\xi]$  satisfies 
\[
 \left(1+\delta \ad_\xi\right)\delta d[\xi,\xi] =0.
\]
But if  $\xi$ is small, $\dd F_\xi = 1+\delta \ad_\xi $ is invertible, so $\delta d[\xi,\xi] =0$.

Finally, we turn to the question of the miniversal family. For simplicity, we discuss only the unobstructed case, i.e, the case
when  $\cK=\cH^1$. If we fix a basis of $\cH^1$, as above, then 
$\widehat{\cO}_{(\cK,0)} \simeq  \CC\llbracket t_1,\ldots, t_d \rrbracket$.
The inverse of the Kuranishi map gives a formal family of deformations of $\Def_L(\CC)$ over $(\cH^1,0)$,
\[
 \Gamma\in L^1\widehat{\otimes} \CC\llbracket t_1,\ldots, t_d \rrbracket, 
\]
\[\Gamma = \sum _{k=1} \Gamma_k,\  \Gamma_k := \sum_{ |J|=k}t^J \Gamma_J,\]
$J$ is a multi-index,  and $\Gamma_k$ are determined inductively: for $x= \sum_{i}t_i\eta_i\in\cH^1$,
\begin{equation}\label{series}
 \Gamma_1(x) = x,\ \Gamma_2(x) = -\frac{1}{2}\delta [\Gamma_1,\Gamma_1],\ \ldots, \Gamma_k = -\frac{1}{2}\delta\sum_{n=1}^{k-1}[\Gamma_n,\Gamma_{k-n}].
\end{equation}
This series has been known for a long time (in various contexts and different levels of generality), 
see e.g.,  \cite{kuranishi_complete}, \cite{kodaira_defo},  \cite{stasheff}, \cite{fukaya_def}, and can be thought of as a reincarnation of
Picard's method of solving ODE's by  iterations.

In the case of a normed dgla, one  shows first   that the series converges (in $\widehat{L}$) in a sufficiently small poly-disk around the origin.
Next,   using elliptic estimates, one proves  that the family can be modified so that the convergence takes place  in $L$.
The prototypical example is the  Kodaira-Spencer dgla,  and the convergence was proved in
\cite{kodaira_nirenberg_spencer}.
In  \cite{itagaki} the author proves  the convergence of this series in the  case of the Barannikov-Kontsevich construction, which contains our setup as
a special case.

\vfill
\eject
\section{Glossary of Notation}\label{conventions}
\noindent $\Art$: the category of local Artin $\CC$-algebras with residue field $\CC$\\
$A$: an Arting ring\\
$\cA^{p}$, $\cA^{p,q}$: sheaves of smooth forms of type $p$ (resp. $(p,q)$) \\
$A^{p}=H^0(X,\cA^p)$, $A^{p,q}=H^0(X,\cA^{p,q})$: global sections\\
$\ad_u=[u,\ ]=\ad u$\\
$\ad\bP = \bP\times_{\ad}\fg$: the adjoint bundle of $\bP$\\

\noindent
$\cB_\fg = H^0(X,\ft \otimes K_X/W)\simeq \oplus_i H^0(X,K_X^{m_i+1})$ the Hitchin base \label{hitchin_base}\\
$\cB^i\subset L^i$ boundaries; used only in Appendix \ref{Kuranishi_Theory}\\
$B_\epsilon$: Ball of radius $\epsilon$, Section \ref{symplectic_kur}\\

\noindent
$\cC^i\subset L^i$: a complement to $\ker d_i$, used only in Appendix \ref{Kuranishi_Theory}\\

\noindent
$d$ or $d_i$: differentials of a complex (always increasing the degree)\\
$\dd$: differential of a map\\
$\Def_{L^\bullet}$: the deformation functor of a dgla $L^\bullet$\\
$\delta$: splitting of a dgla, Appendix \ref{Kuranishi_Theory} \\
$\Delta^+$: positive simple roots, Section \ref{Lie}\\

\noindent
$e_i$: ``upper nilpotent'' Chevalley generators $\{e_i,h_i,f_i\}$\\
$e_{m_i}$: basis vectors for the 1-dimensional subspaces $\fg_{m_i,i}$ \\
$\mend$ (resp. $\send$): Endomorphisms (resp. sheaf endomorphisms);\\
$\mend_m$: $m$-th graded piece of $\mend$\\

\noindent
$\fart$: functors $\cF: \Art\to \sets$ for which $\cF(\CC)=\{\ast\}$\\
$f_i$: ``lower nilpotent'' Chevallety  generators $\{e_i,h_i,f_i\}$\\
$\Phi$:   trivialisation of the symplectic Kuranishi slice \\

\noindent
$G$: simple complex Lie group \\
$\bG$: Green's operator\\
$\fg=\textrm{Lie }G$: simple Lie algebra\\
$\fg_o=Im(\ad_x)\cap Im(\ad_y)$\\
$\fg_m$: $m$-th graded piece of $\fg$ with respect to the principal grading\\
$\fg_{k,i}=\fg_k\cap W_{m_i}$\\
$\fg_x=\fz(x)$: the centraliser of $x\in\fg$\\
$g$: a Riemannian metric on the curve $X$\\
$g_{\overline{\CC}}$: the anti-linear extension of $g$ to $T_{X,\CC}$\\
$g_X$: the genus of the curve $X$\\
$\Gamma = \Phi^{-1}$: the (formal) inverse of the trivialisation of the symplectic Kuranishi slice\\
$\Gamma$: the global section functor\\

\noindent
$\bH$ harmonic projection, $\bH'$, $\bH''$ the two components of $\bH$\\
$\cH$ harmonic representatives of cohomology\\
$\fh$: one of the elements of a principal $\slt$-subalgebra $\{y,\fh,x\}$\\
$\mathpzc{h}$: Coxeter number (largest exponent) of $G$\\
$h$: Hermitian metric on $T_X^\vee$\\
$\bh$: the matrix of the Hermitian metric $h$, a positive real-valued function\\
$(h,v)\in L'\oplus L''=L^1$ a typical element \\
$h_i$: semisimple elements among the Chevalley generators $\{e_i,h_i,f_i \}$ \\

\noindent 
$I$ bundle of centralisers\\
$\io: \Gamma(\send_m(\ad\bP)\otimes K_X^{-m})\hookr \mhom_{C^\infty} (A^{\bullet}(\ad\bP),A^{\bullet}(\ad\bP\otimes K_X^{-m}))$\\
$\io', \io''$: the canonical inclusions of  $L'$ and $L''$ into $L=L'\oplus L''$, Section \ref{symplectic_kur}\\

\noindent 
$\bk:\ft/W\to \fg$: Kostant section, Section \ref{centralisers}\\
$\KK$: formal Kuranishi functor, Appendix \ref{Kuranishi_Theory}, Section 	\ref{symplectic_kur}\\
$\cK$: analytic Kuranishi functor\\
$K_X$ canonical bundle of $X$\\
$\kappa$:  Killing form\\

\noindent 
$l=\textrm{rk}(\fg)$:  the rank of $\fg$\\

\noindent
$m_i$: the exponents of $\fg$\\
$M=\cQ^{-1}(0)$: the vanishing set of the  quadric $\cQ$, Section   \ref{introduction}.\\
$M_{Dol}(G)$: the  Dolbeault moduli space \\
 $\fm_A$: the maximal ideal of $A$\\
$\MC(L)=\cQ^{-1}(0)$: Maurer-Cartan elements of  a dgla $L$,  Section \ref{prelim}\\
$\MC_L$,  the Maurer-Cartan  functors of $L$, Section \ref{prelim}; $\MC_L(A)=\MC(L\otimes A)$\\

\noindent
$o$: marked point\\
$\cO_X$: the structure sheaf of $X$\\

\noindent
$\{p_1,\ldots, p_l\}$: basis of homogeneous invariant polynomials on $\fg$\\
$P$: splitting of $\ad_y$  determined by the choice of principal $\slt$\\
$\bP$ the uniformising Higgs bundle, $\bP=\bF\times_{\ad(\varrho)}G$; also a principal bundle (in general)\\
$pr$: the canonical projection $\MC_L\to \Def_L$\\
$pr^k_n$: the projection $\fg_k\to \fg_{k,n}$, associated with a choice of $\slt$-subalgebra, Section \ref{Lie}\\
$\pi$: a projection $L^2\to \img d'$ in a dgla with a splitting as in Section \ref{symplectic_kur}.\\
$\pi',\pi''$: projections to the two factors $L=L'\oplus L''$ in a dgla with a decomposition as in Section \ref{symplectic_kur}.\\
$pr_i$: projection onto  the $i$-th  factor in a Cartesian product\\

\noindent
$Q$: splitting of $\ad_x$ determined by the choice of principal $\slt$\\
$\cQ$: a quadric; also the Maurer-Cartan quadric $\cQ(u)=du +\frac{1}{2}[u,u]$\\

\noindent 
$\varrho:SL(2,\CC)\to G$ principal embedding\\
$\rho$: Weyl vector, $\rho^\vee$: dual Weyl vector (half the sum of positive coroots)\\

\noindent 
$\SSS$: formal symplectic Kuranishi slice, Section 	\ref{symplectic_kur}\\
$\cS$: analytic symplectic Kuranishi slice\\
$s$: shift, Section \ref{main_example}\\
$\bs: \fz(y)\hookr \fg$: affine-linear map, a variant of Kostant's section\\
$\fs:\cB_\fg\to M_{Dol}(G)$ Hitchin's section\\
$\Sigma\subset \fg$: Kostant's slice,  Section \ref{Lie}\\

\noindent
 $\ft\subset \fg$ Cartan subalgebra

\noindent 
$\underline{V}$: the functor $\underline{V}(A)=V\otimes \fm_A$,   $V$  a vector space \\

\noindent 
$\omega_{can}$: the canonical symplectic form on  $V\times W$, where $V$ and $W$ are  two spaces  in (weak) duality  \\
$\omega$ or $\omega_X$: K\"ahler form on $X$\\
$W$: Weyl group \\
$W_i\subset \fg$: irredicible representations for the principal $\slt$ action on  $\fg$\\
$\cW_\bullet\fg$ : Deligne filtration on $\fg$,  Section \ref{Lie}\\

\noindent 
$x\in \fg$: regular nilpotent, part of a principal $\slt$-subalgebra $\{x,h,y\}$\\
$X$ smooth projective curve (over $\CC$) of genus at least two\\

\noindent 
$y\in \fg$: regular nilpotent, part of a principal $\slt$-subalgebra $\{x,h,y\}$, $y=\sum_i f_i$\\
$\YY$: formal  Kuranishi slice  Appendix \ref{Kuranishi_Theory}, Section 	\ref{symplectic_kur}\\
$\cY$: analytic Kuranishi slice,  Appendix \ref{Kuranishi_Theory}, Section \ref{symplectic_kur}\\

\noindent 
$\fz$: centraliser\\
 $\fz(x)=\fg_x=\oplus \fg_{m_i,i} $\\
 $\fz(y)=\fg_y= \oplus \fg_{-m_i,i}$ \\
 $\zeta$: theta-characteristic\\

\noindent 
 $1_m$: the canonical section of $\cO_X\simeq K^m_X\otimes K^{-m}_X$

\bibliographystyle{alpha}
\bibliography{biblio}

\end{document}